\documentstyle[11pt]{article}
\renewcommand{\baselinestretch}{1.3}

\newtheorem {th}{Theorem}[section]
\newtheorem {lem}[th]{Lemma}
\newtheorem {pr}[th]{Proposition}
\newtheorem {cor}[th]{Corollary}
\newtheorem{defn}{Definition}

\newtheorem{conj}{Conjecture}
\def\Cox{\hfill \Box}

\def\deq{\, {\stackrel {def} {=}}}

\def\ul{\underline}
\def\sf{\sigma\mbox{-field}}
\def\ee{\epsilon}

\def\iid{\mbox{i.i.d.}}

\def\E{{\bf{E}}}
\def\P{{\bf{P}}}

\def\R{{\bf{R}}}
\def\Z{{\bf{Z}}}
\def\S{{\cal{S}}}

\def\F{{\cal{F}}}

\def\|{\, | \, }
\def\one{{\bf 1}}

\def\ct{{\tilde{c}}}

\def\XT{{\tilde{X}}}
\def\ST{{\tilde{S}}}
\def\notlr{\not\leftrightarrow}
\def\com{\mbox{ and } \rho \leftrightarrow}
\def\Cap{\mbox{Cap}}
\def\xx{{\bf x}}
\def\yy{{\bf y}}
\def\zz{{\bf z}}
\def\ww{{\bf w}}

\def\po{\unlhd}
\def\ch{{c_{12}}}
\def\simp{\bigtriangleup}

\begin{document}
 
\setcounter{equation}{0}

\begin{titlepage}
\begin{center}
{\large \bf CRITICAL RANDOM WALK IN RANDOM ENVIRONMENT \\ ON TREES} \\
\end{center}
\vspace{1ex}
\begin{center}
{\sc Robin Pemantle \footnote{Research supported in part by a
National Science Foundation postdoctoral fellowship} $\! ^, \!$ \footnote{Department of
Mathematics, University of Wisconsin-Madison, Van Vleck Hall, 480 Lincoln
Drive, Madison, WI 53706}
 and 
Yuval Peres} \footnote{Current address: Department of Statistics,
University of California, Berkeley, CA 94720} \\
 University of Wisconsin-Madison and Yale University
\end{center}

\vfill

\centerline{\bf Abstract} 
\noindent We study the behavior of Random Walk in Random Environment (RWRE)
on trees in the critical case left open in previous work.  
Representing the random walk by an electrical network, we 
 assume that the ratios of 
resistances of neighboring edges of a tree $\Gamma$ are i.i.d.
random variables whose logarithms have mean zero and finite
variance.  Then the resulting RWRE is transient if simple random
walk on $\Gamma$ is transient, but not vice versa. We obtain 
general transience criteria for such walks, which are sharp for
 symmetric trees of polynomial growth.  In order
to prove  these criteria, we establish
 results on boundary crossing by tree-indexed random walks.  These results
rely on comparison inequalities for percolation processes on trees
and on some new estimates of boundary crossing probabilities for 
ordinary mean-zero finite variance random walks in one dimension,
which are of independent interest.  
\vfill

\noindent{Keywords:} tree, random walk, random environment,
  random electrical network, tree-indexed process, percolation,
 boundary crossing, capacity, Hausdorff dimension.

\noindent{Subject classification: } Primary: 60J15.  Secondary: 60G60, 60G70, 
60E07.

\end{titlepage}

\section{Introduction}

Precise criteria are known for the transience of simple random walk on a
tree (in this paper, a tree is an infinite, locally finite,
rooted acyclic graph and has no leaves, i.e. no vertices of degree one).
See for example Woess (1986), Lyons (1990) or Benjamini and Peres (1992a).
How is the type of the random walk affected if the transition probabilities
are randomly perturbed?  Qualitatively we can say that if this perturbation
has no ``backward push'' (defined below), then the random walk tends to
become more transient; the primary aim of this paper is to establish a
quantitative version of this assertion.

  Designate a vertex $\rho$ of
the tree as its root.  For any vertex $\sigma \neq \rho$, denote by 
$\sigma'$ the unique neighbor of $\sigma$ closer to $\rho$ ($\sigma'$ is
also called the parent of $\sigma$).  An {\em environment\/} for random walk on
a fixed tree, $\Gamma$, is a choice of transition probabilities 
$q(\sigma, \tau)$ on the vertices of $\Gamma$, with $q(\sigma , \tau) > 0$
if and only if $\sigma$ and $\tau$ are neighbors.  When these transition
probabilities are taken as random variables, the resulting mixture of
Markov chains is called {\em Random Walk in Random Environment} (RWRE).
Following Lyons and Pemantle (1992), we study
random environments under the homogeneity condition 
\begin{equation} \label{eq1}
\mbox{The variables } X (\sigma) = \log \left ( {q(\sigma' , \sigma) 
   \over q(\sigma' , \sigma'')} \right ) \mbox{ are i.i.d. for } |\sigma|
   \geq 2 ,
\end{equation}
where $|\sigma|$ denotes the distance from $\sigma$ to $\rho$.  Let
$X$ denote a random variable with this common distribution. 

There are several motivations for studying RWRE under the condition (\ref{eq1}):
\begin{itemize} 
\item For nearest-neighbour RWRE on the integers, the assumption
   (\ref{eq1}) is equivalent to assuming that the transition probabilities
   themselves are i.i.d. .  The first
result on RWRE was obtained by Solomon (1975), who showed that when
$X(\sigma)$ have mean zero and finite variance, then RWRE on the integer
line is recurrent, while it is transient if $E(X)>0$.
 Thus in determining the type of the RWRE, $X$ plays the primary role.
The integer line is the simplest infinite tree, and 
Theorem 2.1 below determines almost exactly the class of trees for which the
same criterion applies. An assumption that  random variables 
analogous to $X(\sigma)$
are stationary (and in particular, identically distributed) is also crucial
in the work of Durrett (1986), which extended the RWRE results of Sinai (1982)
to the multidimensional integer lattice.

\item In terms of the associated resistor network, (\ref{eq1}) means that 
the ratios of the resistances of adjacent edges in $\Gamma$ are i.i.d;
such networks are useful for determining the Hausdorff measures
of certain random fractals-- See Falconer (1988) and Lyons (1990).
The logarithms of the resistances in such a network form
a tree-indexed random walk. (A precise definition of such walks is given below.)
This  structure  appears in a variety of settings: \newline
As a generalization of branching random walk (Joffe and Moncayo 1973);
in a model for ``random distribution functions'' (Dubins and Freedman 1967);
in the analysis of game trees (Nau 1983); 
in studies of random polymers (Derrida and Spohn 1988) and in 
first-passage percolation (Lyons and Pemantle (1992), 
  Benjamini and Peres (1994b)).
\item The tools developed to analyse RWRE satisfying the assumption 
 (\ref{eq1}) are also useful when that assumption is relaxed, e.g. to allow
for some dependence between vertices which are ``siblings''.
In Section 7 we describe an application to certain
  {\em reinforced random walks} which may be reduced to a RWRE.
\end{itemize}

  The main 
result of Lyons and Pemantle (1992) is that RWRE on $\Gamma$ is
a.s. transient if the {\em Hausdorff dimension}, denoted $dim(\Gamma)$, is 
strictly greater than the {\em backward push} 
\begin{equation} \label{eq2}
\beta (X) \deq - \log \min_{0 \leq \lambda \leq 1} \E e^{\lambda X} 
\end{equation}
and a.s. recurrent if $dim(\Gamma) < \beta (X)$.  (The definition of
$dim(\Gamma)$ will be given  in Section~2; the quantity
$e^{dim(\Gamma)}$ is called the branching number of $\Gamma$ in
papers of R. Lyons; the backward push $\beta$ is zero whenever
$X$ has mean zero.) 

While subsuming previous results in Lyons (1990)
and Pemantle (1988), these criteria leave some interesting cases
unresolved.  For instance if $\E X = 0$ (the random environment
is ``fair'') then one easily sees that $\beta (X) = 0$, so the above
criteria yield transience of the RWRE only when $\Gamma$ has positive
Hausdorff dimension, which in particular implies exponential growth.
In fact, much smaller trees suffice for transience of the RWRE in this
case, at least if $X$ has a finite second moment (Theorem~\ref{th2.1}
below).  In particular, this RWRE is a.s. transient whenever simple
random walk on the same tree is transient.  To illustrate the difference
between old criteria such as exponential growth and the criteria
set forth in this paper, we limit the discussion for the rest of the
introduction to {\em spherically symmetric trees}, i.e. trees determined by
a growth function $f : \Z^+ \rightarrow \Z^+$ for which every vertex
at distance $n$ from the root has degree $1 + f(n)$.  Note, however, that
much of the interest in these results stems from their applicability
to nonsymmetric trees; criteria for general trees involve the notion
of capacity and are deferred to the next section.

Assume that $\Gamma$ is spherically symmetric and that the variables
$X(\sigma)$ in~(\ref{eq1}) have mean zero and finite variance
(the assumption of finite variance is plausible since the $X(\sigma)$
are logs of ratios, so the ratios themselves may still have large tails).  
Our first result, Theorem~\ref{th2.1}, is that the RWRE is 
almost surely transient if
\begin{equation} \label{eq3}
\sum_n n^{-1/2} |\Gamma_n|^{-1} < \infty
\end{equation}
and this condition is also necessary, provided that the regularity condition
\begin{equation} \label{eq4}
\sum_n n^{-3/2} \log |\Gamma_n| < \infty
\end{equation}
holds, where $|\Gamma_n|$ is the cardinality of the $n^{th}$ level of $\Gamma$.
We conjecture that condition~(\ref{eq3}) is necessary as well as
sufficient for transience of RWRE.
To see why this is a natural conjecture, and to point out that the 
randomness makes the walk more transient, compare this to the 
following known result: simple random walk on a spherically
symmetric tree is transient if and only if $\sum |\Gamma_n|^{-1}
< \infty$.  

The key to proving the above result is an analysis of {\em tree-indexed
random walks}.  These are random fields $\{ S (\sigma) : \sigma
\in \Gamma \}$ defined from a collection of $\iid$ real random
variables $\{ X (\sigma) : \sigma \in \Gamma, \sigma \neq \rho \}$ by 
letting $S(\sigma)$ be the sum of $X(\tau)$ over vertices $\tau$ on
the path from the root to $\sigma$.  Note that when $\Gamma$ is
a single infinite ray (identified with the positive integers) then
this is just an ordinary random walk; when $\Gamma$ is the family tree
of a Galton-Watson branching process, this is a branching random walk.
Random walks indexed by general trees first appeared in Joffe and Moncayo
(1973).  The motivating question for the study of tree-indexed random walks
(cf. Benjamini and Peres (1994a, 1994b)) is this: when is $\Gamma$ large enough
so that for a $\Gamma$-indexed random walk, the values of $S(\sigma)$
along at least one ray of $\Gamma$ exhibit a prescribed behavior atypical for 
an ordinary random walk?  

In this paper we prove several results in this direction, one of which
we now describe, and apply them to RWRE on trees.  
In the special case where the variables $X(\sigma)$ take only the
values $\pm 1$ with equal probability, Benjamini and Peres (1994b)
obtained conditions for the existence of a ray in $\Gamma$ along which
the partial sums $S(\sigma)$ tend to infinity.  In particular, for
spherically symmetric trees,~(\ref{eq3}) suffices for the existence
of such a ray while the condition
\begin{equation} \label{insert}
\liminf_{n \rightarrow \infty} n^{-1/2} |\Gamma_n| > 0
\end{equation}
is necessary.  In Theorem~\ref{th2.2}
below, this result is extended to variables $X(\sigma)$ with zero mean and
finite variance, and also sharpened.  For spherically symmetric
trees, we show that~(\ref{eq3}) is necessary and sufficient for the existence
of a ray along which $S(\sigma) \rightarrow \infty$.  

As is well known, transience of a reversible Markov chain is equivalent to
finite resistance of the associated resistor network, where the
transition probabilities from any vertex are proportional to the
conductances (reciprocal resistances); see for example Doyle and Snell
(1984).  For an environment satisfying~(\ref{eq1}), the 
conductance attached to the edge between $\sigma'$ and $\sigma$ is
$e^{S(\sigma)}$, where $\{ S(\sigma) \}$ is the $\Gamma$-indexed random
walk with increments $\{ X(\tau) : \tau \neq \rho \}$.  Since finite 
resistance is a tail event, transience of the environment satisfies a
zero-one law.  In particular, the network will have finite resistance
whenever a ray exists along which $S(\sigma) \rightarrow \infty$
sufficiently fast so that $e^{-S(\sigma)}$ is summable.  In this way 
Theorem~\ref{th2.2} yields Theorem~\ref{th2.1}.  

For completeness, we state here a result from Pemantle (1992) about
the case where the $\iid$ random variables $\{ X (\sigma) \}$
have negative mean and the
backward push $\beta (X)$ is positive.  For a spherically 
symmetric tree $\Gamma$, the result of Lyons and Pemantle (1992) 
yields recurrence of the RWRE if
$$\liminf_{n \rightarrow \infty} e^{-n\beta} |\Gamma_n| = 0$$
and transience if
$$|\Gamma_n| \geq Ce^{n(\beta + \ee)}$$
for some $C,\ee > 0$ and all $n$.  Here, analyzing the critical
case is more difficult, but assuming a regularity condition on the random
environment it can be shown that the boundary between
transience and recurrence occurs when 
$$ |\Gamma_n| \approx e^{\beta n + c n^{1/3}} .$$
Here, unlike in the mean zero case, randomness makes the RWRE
more recurrent, since the known necessary and sufficient
condition for transience of RWRE when $X (\sigma) = - \beta$ a.s. is that
$$\sum e^{n \beta} |\Gamma_n|^{-1} < \infty .$$

The rest of the paper is organized as follows.  Precise statements
of our main results are collected in the next section.  Some estimates
for ordinary, mean zero, finite variance random walks that will
be needed in the sequel are collected in Section~3.  Some of these,
along the lines of Woodroofe (1976), may be of independent interest.
In particular, we determine the rate of growth of a mean-zero, finite
variance random walk conditioned to remain positive; this sharpens
considerably a result of Ritter (1981).  
  Section~4 explains
the second moment method for trees, developed by R. Lyons.  Here we
have some new results comparing dependent and independent
percolation and comparing spherically symmetric trees to nonsymmetric
trees of the same size (Theorem~\ref{th4.2}).
After these preliminaries, tree-indexed random walks are discussed
in Section~5, along with an example in which the increments are 
symmetric stable random variables.  The example shows that the {\em
sustainable speed} of a $\Gamma$-indexed random walk
with a given increment distribution may not be determined
by the dimension of $\Gamma$.  Also, this answers a question of R. Lyons
(personal communication) by providing an RWRE satisfying~(\ref{eq1}) which 
is transient on a tree of polynomial growth, even though 
$\E X(\sigma) < 0$.  The RWRE with no backward push is discussed in 
Section~6 and an application to reinforced random walk is described
in Section~7.  

\section{Statements of results}

\setcounter{equation}{0}

We begin with some definitions.  Recall that all our trees are 
infinite, locally finite, rooted at some vertex $\rho$,
and have no leaves.  We use the
notation $\sigma \in \Gamma$ to mean that $\sigma$ is a vertex of
$\Gamma$.
\begin{quote}
1.  An infinite path from the root of a tree $\Gamma$ is called a 
{\em ray} of $\Gamma$.  We refer to the collection of all rays as
the {\em boundary}, $\partial \Gamma$, of $\Gamma$.

2.  If a vertex $\tau$ of $\Gamma$ is on the path connecting the
root, $\rho$, to a vertex $\sigma$, then we write $\tau \leq \sigma$.
For any two vertices $\sigma$ and $\tau$, let $\sigma \wedge \tau$
denote their greatest lower bound, i.e. the vertex where the paths
{}from $\rho$ to $\sigma$ and $\tau$ diverge.  Similarly, the vertex
at which two rays $\xi$ and $\eta$ diverge is denoted $\xi \wedge \eta$.

3.  A set of vertices $\Pi$ of $\Gamma$ which intersects every ray of
$\Gamma$ is called a {\em cutset}.  

4.  Let $\phi : \Z^+ \rightarrow \R$ be a decreasing positive function 
with $\phi (n) \rightarrow 0$ as $n \rightarrow \infty$.  The
Hausdorff measure of $\Gamma$ in gauge $\phi$ is 
$$\liminf_{\Pi} \; \sum_{\sigma \in \Pi} \phi (|\sigma|) ,$$
where the liminf is taken over $\Pi$ such that the distance from $\rho$ to 
the nearest vertex in $\Pi$ goes to infinity.  The supremum over $\alpha$
for which $\Gamma$ has positive Hausdorff measure in gauge $\phi (n)
= e^{-n\alpha}$ is called the {\em Hausdorff dimension} of $\Gamma$.
Strictly speaking, this is the Hausdorff dimension of the boundary of
$\Gamma$ in the metric $d(\xi , \eta) = e^{-|\eta \wedge \xi|}$.  
For spherically symmetric trees, this is just the liminf exponential
growth rate; for general trees it may be smaller.

5.  Hausdorff measure may be defined for Borel subsets $A \subseteq 
\partial \Gamma$ by only requiring the cutsets $\Pi$ to intersect 
all rays in $A$.  Say that $\Gamma$ has $\sigma$-finite Hausdorff
measure in gauge $\phi$ if $\partial \Gamma$ is the union of 
countably many subsets with finite Hausdorff measure in gauge $\phi$.

6.  Say that $\Gamma$ has {\em positive capacity} in gauge $\phi$ if there
is a probability measure $\mu$ on $\partial \Gamma$ for which the
{\em energy} 
$$I_\phi (\mu ) = \int_{\partial \Gamma} \int_{\partial \Gamma}
    \phi (|\xi \wedge \eta|)^{-1} \; d\mu (\xi) \, d\mu (\eta) $$
is finite.  The infimum over probability measures $\mu$ of this
energy is denoted by $1 / \Cap_\phi (\Gamma)$.  
\end{quote}

An important fact about capacity and Hausdorff measure, proved
by Frostman in 1935, is that
$\sigma$-finite Hausdorff measure in gauge $\phi$ implies zero
capacity in gauge $\phi$; the converse just barely fails; 
c.f.\ Carleson (1967 Theorem~4.1).
This gap is either the motivation or the bane of
much of the present work, since many of our criteria would be necessary
and sufficient if zero capacity were identical to $\sigma$-finite
Hausdorff measure.

\begin{th}[proved in Section 6] \label{th2.1}
Suppose that $\iid$ random variables $\{ X(\sigma) : 
\rho \neq \sigma \in \Gamma \}$
are used to define an environment on a tree $\Gamma$ via~(\ref{eq1}),
i.e. the edge from $\sigma'$ to $\sigma$ is assigned the conductance
$$\prod_{\rho < \tau \leq \sigma} e^{X(\tau)} .$$
Assume that $X(\sigma)$ have zero mean and finite variance.
\begin{quote}
$(i)$  If $\Gamma$ has positive capacity in gauge $\phi(n) = n^{-1/2}$,
then the resulting RWRE is transient.

$(ii)$  If $\Gamma$ has zero Hausdorff measure in the
same gauge, then the RWRE is recurrent.

$(iii)$  If $\Gamma$ satisfies the regularity condition
\begin{equation} \label{regularity}
\sum_{n=1}^\infty n^{-3/2} \log |\Gamma_n| < \infty
\end{equation}
then $\sum_{n=1}^\infty n^{-1/2} |\Gamma_n|^{-1} = \infty$ implies
recurrence of the RWRE.  In particular, if $\Gamma$ is spherically symmetric 
and satisfies the regularity condition, then positive capacity in gauge
$n^{-1/2}$ is necessary and sufficient for transience.
\end{quote}
\end{th}

{\bf Remarks:} \\
\noindent{1.}  Lyons (1990) shows that simple random walk ($X(\sigma) =
0$ with probability one) is transient if and only if $\Gamma$ has positive
capacity in gauge $\phi (n) = n^{-1}$.  Thus part $(i)$ of the theorem
justifies the assertion in the introduction that a fair random 
environment makes the random walk more transient.  For spherically
symmetric trees the definitions of Hausdorff measure and capacity
are simpler and the theorem reduces to the conditions 
in~(\ref{eq3})~-~(\ref{insert}).  

\noindent{2.}  Any spherically symmetric tree to which this theorem does
not apply must have zero capacity in gauge $n^{-1/2}$ but fail
the regularity condition; this implies it grows in vigorous
bursts, satisfying $|\Gamma_n| < n^{1/2 + \ee}$ infinitely often, and
$|\Gamma_n| > \exp (n^{1/2 - \ee})$ infinitely often as well. 

\noindent{3.}  If the variables $X(\sigma)$ in~(\ref{eq1}) have 
positive expectation then (trivially) for any tree $\Gamma$ the
RWRE is transient, since the sum of the resistances along any fixed ray
is almost surely finite.  

Part $(i)$ of the theorem is proved by showing that in the mean zero case
there exists a random ray with the same property.  This in turn
is deduced from the next theorem concerning tree-indexed random walks.

\begin{th}[proved in Section 5] \label{th2.2}
Let $\{ X(\sigma) \}$ be $\iid$ random variables indexed by the vertices
of $\Gamma$, and let $S(\sigma) = \sum_{\rho < \tau \leq \sigma} X(\tau)$.
Suppose that $X(\sigma)$ have zero mean and finite variance.  Then
\begin{quote}
$(i)$  If $\Gamma$ has positive capacity in gauge $\phi (n) = n^{-1/2}$
then
$$\P (\exists \xi \in \partial \Gamma \,:\, \forall \sigma \in \xi \;\;
  S(\sigma) \geq 0) > 0 .$$
Furthermore, under the same capacity condition, for every increasing
positive function $f$ satisfying 
\begin{equation} \label{eq5}
\sum_{n=1}^\infty n^{-3/2} f(n) < \infty
\end{equation}
there exists with probability one a ray $\xi$ of $\Gamma$ such that
$S(\sigma) \geq f(|\sigma|)$ for all but finitely many $\sigma \in \xi$.

$(ii)$  If $\Gamma$ has $\sigma$-finite Hausdorff measure in gauge
$\phi (n) = n^{-1/2}$, then 
$$\P (\exists \xi \in \partial \Gamma \,:\, \forall \sigma \in \xi \;\;
  S(\sigma) \geq 0) = 0 .$$
Furthermore, for any increasing $f$ satisfying~(\ref{eq5}), there is
with probability one NO ray $\xi$ such that
$S(\sigma) \geq -f(|\sigma|)$ for all but finitely many $\sigma \in \xi$.

$(iii)$  The conclusions of part $(ii)$ also hold if we assume,
instead of the Hausdorff measure assumption, that $\sum_n n^{-1/2}
|\Gamma_n|^{-1} = \infty$.  
\end{quote}
\end{th}

{\bf Remark:} \\
For spherically symmetric trees, parts $(i)$ and $(iii)$ cover all cases, 
thus proving a sharp dichotomy for tree-indexed random walks: either
there exist rays with $S(\sigma)$ tending to infinity faster than
$n^{1/2 - \ee}$ for all $\ee > 0$, or else along every ray $S(\sigma)$
must dip below $-n^{1/2 - \ee}$ infinitely often.  We believe this
dichotomy holds for all trees but the proof eludes us.  In general,
the condition in part $(iii)$ is not comparable to the condition in
part $(ii)$.  

\section{Estimates for mean zero, finite variance random walk}

\setcounter{equation}{0}

Here we collect estimates for ordinary, one-dimensional, mean zero,
finite variance random walks which are needed in the sequel.  Begin
with a classical estimate whose proof may be found in Feller
(1966), Section XII.8.  

\begin{pr} \label{prFeller}
Let $X_1 , X_2 , \ldots$ be $\iid$, nondegenerate, mean zero, random
variables with finite variance and let $S_n = \sum_{k=1}^n X_k$.
Let $T_0$ denote the hitting time on the negative half-line:
$T_0 = \min \{ n \geq 1 : S_n < 0 \}$.  Then
\begin{equation} \label{eqFeller}
\lim_{n \rightarrow \infty} \sqrt{n} \P (T_0 > n) = c_1 > 0 ,
\end{equation}
and in particular, $c_1' \leq \sqrt{n} \P (T_0 > n) \leq c_1''$
for all $n$.     $\Cox$
\end{pr}

We now determine which boundaries $f(n)$ behave like the horizontal
boundary $f(n) \equiv 0$ in that $\P (S_k > f(k) , k = 1 , \ldots , n)$
is still asymptotically $c n^{-1/2}$.  

\begin{th} \label{th3.1}
With $X_n$ and $S_n$ as in the previous proposition, let $f (n)$
be any increasing positive sequence.  Then 
\begin{quote}
(I) The condition 
\begin{equation} \label{eqsummable}
\sum_{n=1}^\infty n^{-3/2} f(n) < \infty
\end{equation}
is necessary and sufficient for the existence of an integer $n_f$ such that
$$\inf_{n \geq n_f} \sqrt{n} \P (S_k \geq f(k) \mbox{ for } 
   n_f \leq k \leq n) > 0 .$$

(II) The same condition~(\ref{eqsummable}) is necessary and sufficient for
$$\sup_{n \geq 1} \sqrt{n} \P (S_k \geq -f(k) \mbox{ for } 1 \leq k \leq n)
  < \infty .$$
\end{quote}
\end{th}

{\bf Remarks:} \\

\noindent{1.}  Part (I) may be restated as asserting that if $f$
satisfies~(\ref{eqsummable}), then random walk conditioned to stay 
positive to time $n$ stays above $f$ between times $n_f$ and $n$
with probability bounded away from zero as $n \rightarrow \infty$.  This
condition arises in the classical Dvoretzky-Erd\"os test for random walk in
three-space to eventually avoid a sequence of concentric balls of radii
$f(n)$.  Passing to the continuous limit, the absolute value of random walk in
three-space becomes a Bessel (3) process, which is just a
(one-dimensional) Brownian motion conditioned to stay positive.  Proving
Theorem~\ref{th3.1} via this connection (using Skorohod representation,
say) seems more troublesome than the direct proof.

\noindent{2.}  In the case where the positive part of the summands $X_i$
is bounded and the negative part has a moment generating function, 
Theorem~\ref{th3.1} (with further asymptotics) was proved by Novikov
(1981).  He conjectured that these conditions could be weakened.  For
the Brownian case, see Millar (1976).  
Other estimates of this type are given by Woodroofe (1976)
and Roberts (1991).  For their statistical ramifications, see 
Siegmund (1986) and the references therein.  Our estimate can be used 
to calculate the rate of escape of a random walk conditioned to stay
positive {\em forever}; this process has been studied by several
authors -- see Keener (1992) and the references therein.

The proof uses the following three-part lemma.  Let 
$$T_h = \min \{ n \geq 1 : S_n < -h \}$$
denote the hitting time on $(-\infty , - h)$.

\begin{lem} \label{lem ii-iv}
  ~~~
\begin{quote}
(i)  $\P (T_h > n) \leq c_2 h n^{-1/2}$ for all integers $n \geq 1$ and 
real $h \geq 1$.  

(ii)  $\E (S_n^2 \| T_0 > n) \leq c_3 n$ for $n \geq 1$.

(iii)  $\P (T_h > n) \geq c_4 h n^{-1/2}$ for all integers $n \geq 1$ and 
real $h \leq \sqrt{n}$.  
\end{quote}
\end{lem}

{\bf Remarks} (corresponding to the assertions in the lemma):

\noindent{$(i)$}  This estimate is from Kozlov (1976); as we shall see,
it follows immediately from Proposition~\ref{prFeller}.

\noindent{$(ii)$}  In fact we shall verify that
$$\E (S_n^2 \| T_0 > n) \leq 2 n \mbox{ Var}(X_1) + o(n) ,$$
where the constant 2 cannot be reduced in general.  

\noindent{$(iii)$}  Under the additional assumption that $\E |X_1|^3 <
\infty$, this estimate is proved in Lemma~1 of Zhang (1991).  Assuming
only finite variance, as we do, the estimate was known to Kesten
(personal communication) and is implicit in Lawler and Polaski (1992).

{\sc Proof:}  $(i)$  We may assume that $h \leq \sqrt{n}$ and that
$h$ is an integer.  By the central limit theorem there exists an
$r \geq 1$, depending only on the common distribution of the $X_k$,
such that 
$$\P (S_{r^2 h^2} > h) > 1/3 .$$
{}From Proposition~\ref{prFeller} and the FKG inequality
(or the Harris inequality; see Grimmett (1989, Section~2.2)):
$$ \P (S_{r^2 h^2} > h \mbox{ and } T_0 > r^2 h^2) > {c_1'(rh)^{-1} \over 3}
   = c h^{-1} .$$
Consequently
\begin{eqnarray*}
ch^{-1} P (T_h > n) & \leq & P \left [ S_{r^2 h^2} > h \mbox{ and } T_0 > 
   r^2 h^2 \mbox{ and } \sum_{r^2 h^2 + 1}^{r^2 h^2 + k} X_j \geq -h 
   \mbox{ for } 1 \leq k \leq n \right ] \\[2ex]
& \leq & P(T_0 > r^2 h^2 + n) \\[2ex]
& \leq & c_1'' n^{-1/2}
\end{eqnarray*}
which yields the required estimate.

$(ii)$  Consider the minimum of $T_0$ and $n$:
$$ \E (T_0 \wedge n) = \sum_{k=1}^n \P (T_0 \geq k) = \sum_{k=1}^n
   (c_1 + o(1)) k^{-1/2} = 2 (c_1 + o(1)) n^{1/2} .$$
Therefore, using Wald's identity,
$$ \E (S_n^2 \one_{T_0 > n}) \leq \E S_{T_0 \wedge n}^2 = \E X_1^2 
    \E (T_0 \wedge n) = 2 \E X_1^2 (c_1 + o(1)) n^{1/2} .$$
Dividing both sides by $\P (T_0 > n)$ and using Proposition~\ref{prFeller}
gives 
$$\E (S_n^2 \| T_0 > n) \leq (2 + o(1)) n \E X_1^2  .$$
Analyzing the proof, it is easy to see that this is sharp at least when
the increments $X_j$ are bounded.  

$(iii)$  First, we use $(ii)$ to derive the estimate
\begin{equation} \label{eq6}
\E (S_n^2 \| T_h > n) \leq cn 
\end{equation}
with $c$ independent of $h$ and $n$.  By the invariance principle,
$$\inf \{ \P (T_h > n) : n \geq 1 , h \geq \sqrt{n} \} > 0$$
so~(\ref{eq6}) is immediate for $h \geq \sqrt{n}$.  Assume then that
$h < \sqrt{n}$.  Let $A_i$ denote the event that $T_h > n$ and $S_i$
is the last minimal element among $0 = S_0 , S_1 , \ldots , S_n$.  Then
$$\E (S_n^2 \| T_h > n) = \sum_{i=1}^n \P (A_i) \E (S_n^2 \| A_i) .$$
Conditioning further on $S_i$ we see that this is at most
$$ \sup \{ \E (S_n^2 \| A_i , S_i = y) : 1 \leq i \leq n , -h \leq 
   y \leq 0 \} . $$
But by the Markov property,
$$\E (S_n^2 \| A_i , S_i = y) = \E ((y + S_{n-i})^2 \| S_k > 0 \mbox{ for }
    1 \leq k \leq n-i ) .$$
Since $y < 0$ and $y^2 < n$, this gives
$$ \E (S_n^2 \| T_h > n) \leq \sup_{0 \leq i \leq n} 
   \{ n + \E (S_{n-i}^2 \| S_k > 0
   \mbox{ for } 1 \leq k \leq n-i) \} \leq (1+c_3)n $$
with $c_3$ as in $(ii)$, proving~(\ref{eq6}).  

Now letting $A$ denote the event $\{ T_h > n \}$, we have
$$0 = \E S_{T_h \wedge n} = \P (A) \E (S_n \| A) + \P (A^c) \E (S_{T_h}
   \| A^c) $$
and therefore (using~(\ref{eq6}) in the final step):
$${\P (A) \over \P (A^c)} = {\E (-S_{T_h} \| A^c) \over \E (S_n \| A)} 
   \geq {h \over (\E (S_n^2 \| A))^{1/2}} \geq {h \over (cn)^{1/2}} .$$
This establishes $(iii)$.     $\Cox$

{\sc Proof of Theorem}~\ref{th3.1}, part (I): First assume that 
$\sum n^{-3/2} f(n) < \infty$.  Consider the events 
$$V_m = \left \{ S_k > f(2^m) \mbox{ for } 2^{m-1} < k \leq 2^m \right \} .$$
We claim that for any $N \geq 2^{m-1}$,
\begin{equation} \label{eq7}
\P (V_m^c \| T_0 > 4N) \leq  \ct f(2^m) 2^{-m/2} 
\end{equation}
with some constant $\ct > 0$.  Indeed by conditioning on the first $k \in
[2^{m-1} + 1 , 2^m]$ for which $S_k \leq f(2^m)$ one sees that
\begin{eqnarray*}
&& \P (V_m^c \cap \{ T_0 > 4N \}) \\[2ex]
& \leq & \P (T_0 > 2^{m-1}) \max_{k \leq 2^m} \P (S_j - S_k \geq -f(2^m)
   \mbox{ for all } j \in [k+1 , 4N] ) \\[2ex]
& \leq & \P (T_0 > 2^{m-1}) \P (T_{f(2^m)} \geq N) \\[2ex]
& \leq & c_1''c_2 f(2^m) (N 2^{m-1})^{-1/2} .
\end{eqnarray*}
Using Proposition~\ref{prFeller}, this establishes~(\ref{eq7}).  
Since $f$ is nondecreasing, the hypothesis 
$$\sum_n n^{-3/2} f(n) < \infty$$
is equivalent to $\sum_{m=1}^\infty 2^{-m/2} f(2^m) < \infty$.  
Choose an integer $m_f$ such that 
$$ \ct \sum_{m = m_f}^\infty 2^{-m/2} f(2^m) < {1 \over 2} .$$
By~(\ref{eq7}), for every $M > m_f$ we have 
$$\sum_{m = m_f}^M \P (V_m^c \| T_0 > 2^{M+1}) \leq {1 \over 2}$$
and hence
$$\P \left ( \bigcap_{m = m_f}^M V_m \| T_0 > 2^{M+1} \right ) \geq 
    {1 \over 2} .$$
Taking $n_f = 2^{m_f}$ and recalling the definition of $V_m$ concludes the
proof.

For the converse, first recall a known fact about mean zero, finite
variance random walks, namely that
\begin{equation} \label{eq8}
\lim_{R \rightarrow \infty} \sup_{h \geq 0} \P (h + S_{T_h} \leq -R) = 0 .
\end{equation}
In other words the amounts by which $\{ S_n \}$ {\em overshoots} the
boundary $-h$ are tight as $h$ varies over $(0,\infty)$.  Indeed
the overshoots for the random walk $\{ S_n \}$ are the same as for the
associated renewal process $\{ L_n \}$ of descending ladder random
variables, where $L_1$ is the first negative value among $S_1 , S_2 ,
\ldots$, and in general, $L_{n+1}$ is the first among $\{ S_k \}$
which is less than $L_n$.  The differences $L_{n+1} - L_n$ are
$\iid$; they have finite first moment if and only if $\E X_1^2 < \infty$
(see Feller (1966), Section XVIII.5).  In this case, the overshoots are
tight since by the renewal theorem, they converge in distribution
(see Feller (1966), Section XI.3).  This proves~(\ref{eq8}).

Some new notation will be useful:
\begin{defn}
For any function $f(n)$ and any random walk $\{ S_n \}$, let 
$A(f;a,b)$ denote the event that $S_n \geq f(n)$ for all
$n \in [a,b]$.  Let $A(f;b)$ denote $A(f;1,b)$.
\end{defn}
Proceeding now to the proof itself, it is required to 
prove~(\ref{eqsummable}) from the assumption that for some $n_f \geq 1$,
\begin{equation} \label{eq9}
\inf_{n \geq n_f} \P ( A(f;n_f,n) \| T_0 > n) > 0 .
\end{equation}
It may be assumed without loss of generality that $f(n) \rightarrow \infty$,
since otherwise there is nothing to prove; also, by changing
$f$ at finitely many integers, it may be assumed, without affecting the
condition~(\ref{eqsummable}) we are trying to prove, that~(\ref{eq9})
holds with $n_f = 1$.  

Impose the restriction $f(n) \leq \sqrt{n}$; this restriction will
be removed at the end of the proof.  Let $c_6$ be the infimum
of probabilities $\P (A(f;n) \| T_0 > n)$, which is positive
by~(\ref{eq9}).  The key estimate to proving~(\ref{eqsummable}) is
\begin{equation} \label{eq11}
\P (A(f;2n)^c \| A(f;n) \mbox{ and } T_0 > N) \geq c_7 f(n) n^{-1/2}
\end{equation}
for some $c_7 > 0$, all sufficiently large $n$ and all $N \geq 2n$.
Verifying this estimate involves several steps.
\begin{quote}
\ul{Step 1: Controlling $S_n$, given $A(f;n)$}.  From part~$(ii)$ of
Lemma~\ref{lem ii-iv}, 
$$\E (S_n^2 \| A(f;n)) \leq c_6^{-1} \E (S_n^2 \| T_0 > n) 
   \leq (c_3 / c_6) n .$$  
Therefore
\begin{equation} \label{eq12}
\P (S_n \leq c_8 \sqrt{n} \| A(f;n)) \geq {1 \over 2},
\end{equation}
where $c_8 = 2 c_3 / c_6 > 0$.

\ul{Step 2: Securing a dip below the boundary}.  By the central limit theorem
there exists an integer $n^*$ and a constant $c_9 > 0$ such that
$$\P (S_{2n} - S_n < - c_8 \sqrt{n} ) \geq 4 c_9 \mbox{ for } n \geq n^* .$$
Let $t(n) = \min \{ k > n : S_k < f(n) \}$.  Then nonnegativity
of $f$ and the Markov property imply 
$$\P (t(n) \leq 2n \| A(f;n) , S_n ) \geq 4 c_9 \one_{S_n \leq 
   c_8 \sqrt{n}} ,$$
and hence by~(\ref{eq12}),
\begin{equation} \label{eq13}
\P (t(n) \leq 2n \| A(f;n)) \geq 2 c_9 
\end{equation}
whenever $n \geq n^*$.

\ul{Step 3: Controlling the overshoot}.  Use tightness of the overshoots
to pick an $R > 0$ such that $\P (S_{t(n)} \geq f(n) -
R \| S_n = y) \geq 1 - c_9$
for any $y \geq f(n)$.  Increase $n^*$ if necessary to ensure that
$f(n^*) > 2R$ and hence for all $n \geq n^*$, $\P (S_{t(n)} \geq
f(n) / 2 \| A(f;n)) \geq 1 - c_9$.  Combining this with~(\ref{eq13}) yields
$$\P \left ( t(n) \leq 2n \mbox{ and } S_{t(n)} \geq {f(n) \over 2} 
   \| A(f;n) \right ) \geq c_9 ;$$
thus by Proposition~\ref{prFeller} and the definition of $c_6$ there
is some $c_{10} > 0$ such that for all $n$,
\begin{equation} \label{eq 940102a}
\P (A(f;n) \cap \{ t(n) \leq 2n \} \cap \{ S_{t(n)} \geq {f(n) \over 2}
   \}) \geq c_{10} n^{-1/2} .
\end{equation}

\ul{Step 4: Maintaining positivity}.  From the strong Markov property 
and part~$(iii)$ of Lemma~\ref{lem ii-iv}, the event
$\{ S_k - S_{t(n)} \geq -f(n) / 2 \mbox{ for } k \in [t(n)+1 , t(n) +N] \}$
is independent of the random walk up to time $t(n)$ and has probability at
least $(c_4 / 2) f(n) N^{-1/2}$.  Multiplying by 
the inequality~(\ref{eq 940102a}) proves that
$$\P (A(f;n) \cap A(f;2n)^c \cap \{ T_0 > N \}) \geq 
c_{11} f(n) (nN)^{-1/2} ,$$
where $c_{11} = c_{10} \cdot c_4 / 2$.
Now the key estimate~(\ref{eq11}) follows from Proposition~\ref{prFeller}.
\end{quote}

\noindent{From} here the rest is easy sailing.  If $n > n^*$ and $2^M \geq 2n$ 
then
$$ \P (A(f;2n) \| T_0 > 2^M) \leq (1 - c_7 f(n) n^{-1/2}) \P (A(f;n) 
    \| T_0 > 2^M) . $$
Therefore
$$\P (A(f;2^M) \| T_0 > 2^M) \leq \prod_{\log_2 n^* < m \leq M} 
   1 - c_7 f(2^m) 2^{-m/2} .$$
Recalling that the LHS is bounded away from zero, we
infer that necessarily 
$$\sum_m 2^{-m/2} f(2^m) < \infty .$$  

Having proved summability from~(\ref{eq9}) when $f(n) \leq \sqrt{n}$,
we now remove the restriction.  If the restriction is violated
finitely often, this is easily corrected where it is used in Step~2
by choosing $n^*$ sufficiently large.  If the restriction is violated 
infinitely often, then the above proof works for $g(n) = f(n) \wedge
\sqrt{n}$ to show that 
$$ \sum_{m=1}^\infty 2^{-m/2} \min (f(2^m) , 2^{m/2}) < \infty $$
which contradicts infinitely many violations.  Thus 
$2^{-m/2} f(2^m)$ is summable in any event, which is equivalent
to~(\ref{eqsummable}).    $\Cox$

{\sc Proof of Theorem}~\ref{th3.1}, part (II):
We may assume that $f(n) \uparrow \infty$ since Lemma~\ref{lem ii-iv}
part $(i)$ covers bounded $f$.  Retain the notation $A(f;a,b)$,
{}from the previous proof.  One of the two halves of the equivalence,
$\sup_{n \geq 1} \sqrt{n} \P (S_k \geq -f(k) \mbox{ for } 1 \leq k \leq n)
< \infty$ implies summability of~(\ref{eqsummable}), is easy.
Assume that $\P (A(-f;n)) \leq C n^{-1/2}$ for all 
$n \geq 1$.  Under this assumption, we may repeat the proof 
of~(\ref{eq11}) substituting 0 for the upper boundary, $f$, and
substituting $-f$ for the lower boundary, 0, to yield
\begin{equation} \label{eq18}
\P (A(0;2n)^c \| A(0;n) \cap A(-f;N)) \geq c_7 f(n) n^{-1/2}
\end{equation} 
for all large $n$ and $N \geq 4n$.  Our assumption implies that
the products
$$ \prod_{m=1}^M \P (A(0;2^{m+1}) \| A(0;2^m) \cap A(-f;2^{M+2})) , $$
being greater than $\P (A(0;2^{M+1}) \| A(-f;2^{M+2}))$, must
be bounded below by a positive constant.  From~(\ref{eq18}) it follows
that the product of $1 - c_7 n^{-1/2} f(n)$ is nonzero as $n$
ranges over powers of two, which implies $\sum 2^{-m/2} f(2^m) < \infty$,
completing the proof for this half of the equivalence.   

The other direction rests on an inequality which will
require some work to prove. 
{\sc Claim:} {\em There is a constant $\, \ch \geq 1 \, $ for which}
\begin{equation} \label{new 3.10}
\P \left ( T_0 < 2n \| T_0 \geq n \mbox{ and } A(-f ; N) \right ) 
   \leq \ch f(3n) n^{-1/2}
\end{equation}
{\em provided that\/} $N \geq 4n$. (Observe that when $f(3n)^2 \geq n$ the 
inequality is trivial. Also note that it suffices to establish
 (\ref{new 3.10}) for large $n$, since $\ch$ can be chosen large enough
to render the inequality trivial for small $n$.).

  Assuming~(\ref{new 3.10}) for the moment,
\begin{eqnarray} 
&& \P (T_0 \geq 2^{M+1} \| A(-f; 2^{M+2})) \nonumber \\[2ex]
& \geq &  \P (T_0 \geq 2^{m_0} \| A(-f; 2^{M+2})) \cdot \prod_{m=m_0}^M
   \left ( 1 - \ch f(3 \cdot 2^m) 2^{-m/2} \right ) \nonumber \\[2ex]
& \geq & c_{m_0} \prod_{m=m_0}^M \left ( 1 - \ch f(3 \cdot 2^m) 
   2^{-m/2} \right ) ,  \label{new 3.11} 
\end{eqnarray}
where $m_0$ is large enough so that all the factors on the right
positive.  Since the summability of~(\ref{eqsummable}) is equivalent to
$$ \sum f(3 \cdot 2^m) 2^{-m/2} < \infty , $$
the RHS of~(\ref{new 3.11}) is bounded from below by a positive
constant $c_{13}$ and therefore 
$$\P (A( -f , 2^{M+2})) \leq c^{-1} \P (T_0 \geq 2^{M+1}) .$$
Proposition~\ref{prFeller} then easily yields the inequality
we are seeking (use the smallest $M$ such that $2^M \geq n$).  

Thus it suffices to establish~(\ref{new 3.10}).  This is done by 
cutting and pasting portions of the random walk trajectory. 
To bound the LHS of~(\ref{new 3.10}) we will condition on the time 
$T_0 = k \in [n , 2n)$ of the first negative value of the random walk,
and on the overshoot $S_{T_0} = -y < 0$.  This gives
\begin{equation} \label{new 3.12}
\P ( n \leq T_0 < 2n \mbox{ and } A(-f;N)) \leq P (T_0 \geq n)
   \max_{\begin{array}{c} n \leq k < 2n \\ 0 < y < f(k) \end{array}} 
   \P (A( -f ; k , N) \| S_k = -y) .
\end{equation}
Next observe that, given $S_k = -y$, the events $A(-f; k,N)$ and
$\{ S_{3n} > \sqrt{n} \}$ are both increasing events in the
conditionally independent variables $X_{k+1} , \ldots , X_N$.
Applying the Harris inequality (or FKG) yields
\begin{eqnarray}
&& \P (A(-f; k,N) \mbox{ and } S_{3n} > \sqrt{n} \| S_k = -y) 
   \nonumber \\[2ex]
& \geq & \P (A(-f; k,N) \| S_k = -y) \cdot \P ( S_{3n} > \sqrt{n} 
   \| S_k = -y) \nonumber \\[2ex]
& \geq & c_{14} \P (A(-f; k,n) \| S_k = -y)  , \label{new 3.13}
\end{eqnarray}
where the last inequality uses the fact that $3n-k > n$,
that $0 < y \leq f(3n) < \sqrt{n}$, and the central limit theorem.

Now begins the cutting and pasting.  We shall combine a trajectory
$X_{k+1}^{(1)} , X_{k+2}^{(1)} , \ldots ,  X_N^{(1)}$ in the
event on the LHS of~(\ref{new 3.13}) (called $B_1$ in figure~1)
with two independent random walk trajectories $\{ X_j^{(2)} : 
j \geq 1 \}$ and $\{ X_j^{(3)} : j \geq 1 \}$ depicted in 
figures~2 and~3 respectively.  

Assume without loss of generality that $f(3n)^2$ is an integer.  
Define an event (i.e.\ a subset of sequence space) by
$$B(k,y) = A (y-f; 0, N-k) \cap \{ S_{3n-k} > \sqrt{n} + y \}$$
and observe that $B_1 \cap \{ S_k = -y \}$ may be written
as the intersection of $\{ S_k = -y \}$ 
with the event that the shifted sequence 
$X_{k+1}^{(1)}, X_{k+2}^{(1)}, \ldots$ is in $B(k,y)$.  
Define the
mapping taking the three trajectories $\{ X_j^{(i)} \}$, $i = 1,2,3$
into a trajectory $\{ \XT_j : 1 \leq j \leq N \}$ as follows.
$$\begin{array}{rrcl}
I. & \XT_j = X_j^{(2)} & \mbox{ for } & 1 \leq j \leq f(3n)^2 \\[2ex]
II. & \XT_{f(3n)^2 + j} = X_{k+j}^{(1)} & \mbox{ for } & 1 \leq j 
   \leq 3n - k \\[2ex]
III. & \XT_{f(3n)^2 + 3n - k + j} = X_j^{(3)} & \mbox{ for } & 1 \leq j 
   \leq k - f(3n)^2 \\[2ex]
IV. & \XT_j = X_j^{(1)} & \mbox{ for } & 3n + 1 \leq j \leq N . 
\end{array}$$

\begin{picture}(400,100)
\put(100, 20) {\framebox(200,60){Figures 1 - 4 go here}}
\end{picture}

\noindent We claim that the ``pasted'' trajectory  $\{\XT_j : 1 \leq j \leq N \}$
lies in the event \newline $B_4 \deq A(0;n) \cap A(-f;N)$ depicted in figure~4
whenever the trajectories $\{ X_j^{(1)} \}$, $\{ X_j^{(2)} \}$ and 
$\{ X_j^{(3)} \}$ lie in $B_1$, $B_2$ and $B_3$ respectively
(See figures~1-4 for the definitions of $B_1, B_2$ and $B_3$).
Indeed, let $\ST_j = \XT_1 + \cdots + \XT_j$ and observe that for
$1 \leq j \leq 3n - k$, 
$$ \sum_{i=k+1}^{k+j} X_i^{(1)} \geq - f(k+j) - (-y) \geq - f(3n) .$$
Thus $\ST_{f(3n)^2 + j} \geq \ST_{f(3n)^2} - f(3n) \geq 0$.  This 
verifies that part II of the trajectory in figure~4 satisfies
the requirements to be in $B_4$; the other verifications are immediate.

Since the sequence $\{ \XT_j \}$ is a fixed permutation of 
the three sequences $\{ X_j^{(1)} \}$, $\{ X_j^{(2)} \}$ and 
$\{ X_j^{(3)} \}$, it is still i.i.d.\ and hence
\begin{equation} \label{new 3.15}
\P (B(k,y)) \P (B_2) \P (B_3) = \P (B_1 \| S_k = -y) \P (B_2) \P (B_3) 
   \leq \P (B_4) .
\end{equation}
Now the event $B_2$ is the intersection of two increasing events, so
by the Harris inequality, Proposition~\ref{prFeller} and the central 
limit theorem,
\begin{equation} \label{new 3.16}
\P (B_2) \geq \P (A(0 ; f(3n)^2)) \P (S_{f(3n)^2} \geq f(3n)) \geq
   { c_{15} \over f(3n)}
\end{equation}
for some positive $c_{15}$ when $n$ (and hence $f(3n)$) are sufficiently
large.  

Similarly, by Lemma~\ref{lem ii-iv}~$(iii)$, the CLT and the Harris inequality, 
\begin{equation} \label{new 3.17}
\P (B_3) \geq ({1 \over 2} - o(1)) \P (A(- \sqrt{n} ; k-f(3n)^2 )) 
   \geq c_{16} > 0.
\end{equation}
Together with~(\ref{new 3.15}) and~(\ref{new 3.16}) this yields
\begin{equation} \label{new 3.18}
\P (B_1 \| S_k = -y) \leq c_{17} f(3n) \P (B_4) .
\end{equation}
By~(\ref{new 3.13}), 
$$\P (A( -f; k,N) \| S_k = -y) \leq c_{18} f(3n) \P (B_4) $$
and since this is true for any choice of $y$, we integrate out $y$
to get
$$\P (A( -f; k,N)) \leq c_{18} f(3n) \P (B_4) .$$
Finally, recalling~(\ref{new 3.12}) and the definition of $B_4$, we obtain
\begin{eqnarray*}
\P (n \leq T_0 < 2n \mbox{ and } A(-f; N)) & \leq & c_{18} \P (T_0
   \geq n) f(3n) \P (B_4) \\[2ex]
& \leq & c_{19} n^{-1/2} f(3n) \P (T_0 \geq n \mbox{ and } A(-f ; N))
\end{eqnarray*}
which is equivalent to~(\ref{new 3.10}). 
  This completes the proof of 
Theorem~\ref{th3.1}.   $\Cox$

\section{Target percolation on trees}

\setcounter{equation}{0}

In this section only, it will be convenient to consider finite as well as
infinite trees.  Among finite trees, we allow only those of constant
height, i.e. all maximal paths from the root (still called rays) have the
same length $N$.  Thus the set $\partial \Gamma$ of rays may be identified
with $\Gamma_N$.  The definitions of energy and capacity remain the same.
The definition of Hausdorff measure fails, since the cutset $\Pi$ cannot
go to infinity; replacing the liminf by an infimum defines the {\em 
Hausdorff content}.  

Following Lyons (1992), we consider a very general {\em percolation
process} on $\Gamma$: a random subgraph $W \subseteq \Gamma$ chosen 
{}from some arbitrary distribution on sets of vertices in $\Gamma$.  The
event that $W$ contains the path connecting $\rho$ and $\sigma$ is
denoted $\{ \rho \leftrightarrow \sigma \}$; similarly write $\{ \rho
\leftrightarrow \partial \Gamma \}$ for the event that $W$ contains a
ray of $\Gamma$.  A familiar example from percolation theory is when
each edge $e$ is retained independently with some probability $p(e)$.  
The random component $W$ of this subgraph that contains the root 
is called a {\em Bernoulli percolation} on $\Gamma$.
Another example that is general enough to include nearly
all cases of interest is a {\em target} percolation.  This is defined
{}from a family of $\iid$ real random variables $\{X(\sigma)
: \sigma \neq \rho \}$ by choosing some
closed set $B \subseteq \R^N$ and defining
$$W = \bigcup_{k=0}^N \{ \sigma \in \Gamma_k : (X(\tau_1) , 
   \ldots , X(\tau_k)) \in \pi_k B \} ,$$ 
where $\pi_k$ is the projection of $B$ on the first $k$ coordinates,
the sequence $\rho , \tau_1, \ldots , \tau_k = \sigma$ is the 
path from the root to $\sigma$, and $\rho$ is defined always to be 
in $W$.  The set $B$ is called the target set.
Observe that since $B$ is closed, a ray $\xi = (\rho, \sigma_1 , \sigma_2 , 
\ldots)$ is in $W$ if and only if $(X(\rho) , X(\sigma_1) , \ldots) \in B$.
Letting $\{ X (\sigma) \}$ all be uniform on $[0,1]$ and 
$B = \{ \prod_{j=1}^\infty [0,a_j] \}$ for some $a_j \in [0,1]$
recovers a class of Bernoulli percolations.

The following lemma, which will be sharpened below, is contained in 
the results of Lyons (1992).  Because of notational differences,
the brief proof is included.

\begin{lem} \label{lem4.1}
Consider a percolation in which $\P (\rho \leftrightarrow \sigma) 
= p(|\sigma|)$ for some strictly positive function $p$.  
\begin{quote}
$(i)$  First moment method: $\P (\rho \leftrightarrow \partial \Gamma)$
is bounded above by the Hausdorff content of $\Gamma$ in the gauge
$p(n)$.  If $\Gamma$ is infinite and has zero Hausdorff measure in gauge
$\{ p(n) \}$, then $\P (\rho \leftrightarrow \partial \Gamma) = 0$.

$(ii)$  Second moment method: Suppose further that there is a positive, 
nonincreasing function $g : \Z^+ \rightarrow \R$ such that for any two 
vertices $\sigma , \tau \in \Gamma_n$ with $|\sigma \wedge \tau| = k$,
\begin{equation} \label{eqQB}
\P (\rho \leftrightarrow \sigma \com \tau) \leq {p(n)^2 \over g(k) } .
\end{equation}
Then $\,\P (\rho \leftrightarrow \partial \Gamma) \geq \Cap_g (\Gamma)$.
\end{quote}
\end{lem}

{\bf Remark:} \\
For Bernoulli percolations,~(\ref{eqQB}) holds with equality for $g(k) =
p(k)$.  More generally, when $g(k) = p(k) / M$ for some constant $M > 0$,
the percolation is termed {\em quasi-Bernoulli} (Lyons 1989).  In this
case, 
\begin{equation} \label{eqQB2}
\P (\rho \leftrightarrow \partial \Gamma) \geq {\Cap_p (\Gamma) \over M} .
\end{equation}
Note also that any target percolation satisfies the condition in the
lemma: $\P (\rho \leftrightarrow \sigma) = p(|\sigma|)$.

{\sc Proof:}  $(i)$:  For any cutset $\Pi$,
$$\P (\rho \leftrightarrow \partial \Gamma) \leq \P (\rho \leftrightarrow
  \sigma \mbox{ for some } \sigma \in \Pi) \leq \sum_{\sigma \in \Pi}
  p(|\sigma|) .$$
The assertion follows by taking the infimum over cutsets.  

$(ii)$:  First assume that $\Gamma$ has finite height $N$.  Let $\mu$ be
a probability measure on $\partial \Gamma = \Gamma_N$.  Consider
the random variable 
$$ Y_N = \sum_{\sigma \in \Gamma_N} \mu (\sigma) \one_{\rho \leftrightarrow
  \sigma} .$$
Clearly $\E Y_N = p(N)$ and 
\begin{eqnarray*}
\E Y_N^2 & = & \E \sum_{\sigma , \tau \in \Gamma_N} \mu (\sigma) \mu (\tau)
  \one_{\rho \leftrightarrow \sigma \com \tau} \\[2ex]
& \leq & \sum_{\sigma , \tau \in \Gamma_N} \mu (\sigma) \mu (\tau)
  { p(N)^2 \over g(|\sigma \wedge \tau|)} \\[2ex]
& = & p(N)^2 I_g (\mu) 
\end{eqnarray*}
by the definition of the energy $I_g$.  The Cauchy-Schwartz inequality
gives
$$p(N)^2 = \E (Y_N \one_{Y_n > 0})^2 \leq \E Y_N^2 \P (Y_N > 0)$$
and dividing the previous inequality by $\E Y_N^2$ gives $1 \leq I_g (\mu)
\P (Y_N > 0)$.  Since $\Cap_g (\Gamma)$ is the supremum of $I_g(\mu)^{-1}$
over probability measures $\mu$, it follows that $\P (Y_N > 0) \geq \Cap_g
(\Gamma)$.  The case where $N = \infty$ is obtained form a straightforward
passage to the limit.  $\Cox$

For infinite trees, we are primarily interested in whether $\P (\rho
\leftrightarrow \partial \Gamma)$ is positive.  The above lemma 
fails to give a sharp answer
even for quasi-Bernoulli percolations, since the condition 
$\Cap_p (\Gamma) = 0$ does not imply zero Hausdorff measure.  We believe 
but cannot prove the following.
\begin{conj} \label{conjcapacity}
For any target percolation on an infinite tree with $\P (\rho 
\leftrightarrow \partial \Gamma) > 0$, the tree $\Gamma$ must
have positive capacity in gauge $p(n)$, where $p(|\sigma|) =
\P (\rho \leftrightarrow \sigma)$.
\end{conj}
To see why we restrict to target percolations, let $\Gamma$ be the
infinite binary tree, let $\xi$ be a ray chosen uniformly from the
canonical measure on $\partial \Gamma$ and let $W = \xi$.  Then
$\P (\rho \rightarrow \partial \Gamma) = 1$, but $\Gamma$ has zero
capacity in gauge $p(n) = 2^{-n}$.  Evans (1992) gives a capacity
criterion on the target set $B$ necessary and sufficient for 
$\P (\rho \leftrightarrow \partial \Gamma) > 0$ in the special case
where $\Gamma$ is a homogeneous tree (every vertex has the same degree).
His work was extended by Lyons (1992), who showed that
$\Cap_p (\Gamma) > 0$ was necessary for $\P(\rho \leftrightarrow \partial
\Gamma) > 0$ for all Bernoulli percolations and for 
non-Bernoulli percolations satisfying a certain condition.  

Specific non-Bernoulli target percolations are used in Lyons (1989)
to analyze the Ising model, and in Lyons and Pemantle (1992),
Benjamini and Peres (1994b) and Pemantle and Peres (1994) to determine 
the speed of first-passage percolation.  In the present work, the
special case 
$$B = \{ \xx \in \R^\infty : \sum_{i=1}^n x_i \geq 0 \mbox{ for all } n \}$$
will play a major role.  Unfortunately, this set does not satisfy  
Lyons' (1992) condition.  It is, however, an {\em increasing set},
meaning that if $\xx \geq \yy$ componentwise and $\yy \in B$,  
then $\xx \in B$.  This motivates the next lemma.

\begin{lem}[sharpened first moment method] \label{sharpened}
Consider a target percolation on a tree $\Gamma$ in which the target
set $B$ is an increasing set.  Assume that $p(n) \deq p(\rho
\leftrightarrow \sigma)$ for $\sigma \in \Gamma_n$ goes to zero
as $n \rightarrow \infty$.
\begin{quote}
$(i)$  With probability one, the number of surviving rays (elements of $W$)
is either zero or infinite.

$(ii)$  If $\partial \Gamma$ has $\sigma$-finite Hausdorff measure in the 
gauge $\{ p(n) \}$, then $\P (\rho \leftrightarrow \partial \Gamma) = 0$.
\end{quote}
\end{lem}

\noindent{Remark:}  With further work, we can show that the assumption 
that $B$ is an increasing set may be dropped; since this is the only case 
we need, we impose the assumption to greatly simplify the proof.  The above 
example shows that the ``target'' assumption cannot be dropped.

\noindent{\sc Proof:}  Assume that $\P (\mbox{ finitely many rays survive})
> 0$.  Let $A_k$ denote the event that exactly $k$ rays survive and
fix a $k$ for which $\P (A_k) > 0$.
Let $\F_n$ denote the $\sf$ generated by $\{ X (\sigma) : |\sigma|
\leq n \}$.  Convergence of the martingale $\P (A_k \| \F_n)$ 
shows that for sufficiently large $n$ the probability 
that $\P (\mbox{ exactly $k$ surviving rays} \| \F_n) > .99$ is positive.
Let $\{ x (\sigma) : |\sigma| \leq n \}$ be a set of values for which
$$\P (\mbox{exactly $k$ rays survive} \| X(\sigma) = x (\sigma) :
  |\sigma| \leq n) > 0.99 \, .$$
Totally order the rays of $\Gamma$ in any way; since the probability
of any fixed ray surviving is zero, it follows that there is 
a ray $\xi_0$ such that 
$$\P (\mbox{some $\xi < \xi_0$ survives} \| X(\sigma) = x (\sigma) :
  |\sigma| \leq n) = {1 \over 2} \, .$$
This implies that 
$$\P (\mbox{at least $k$ rays $\xi > \xi_0$ survive} \| 
   X(\sigma) = x (\sigma) : |\sigma| \leq n) \geq 0.49 \, .$$
Since all the $X(\sigma)$ are conditionally 
independent given $X(\sigma) = x (\sigma)$ for $|\sigma| \leq n$, and since
the events of at least one ray less than $\xi_0$ surviving and
at least $k$ rays greater than $\xi_0$ surviving are both
increasing, we may apply the FKG inequality to conclude that
$$\P (\mbox{at least $k+1$ rays survive} \| 
   X(\sigma) = x (\sigma) : |\sigma| \leq n) \geq {0.49 \over 2} \, .$$
This contradicts the choice of $x(\sigma)$,
so we conclude that 
$$\P (\mbox{ finitely many rays survive}) = 0.$$

For part~$(ii)$, write $\Gamma = \bigcup_{n=1}^\infty \Gamma_k$
where each $\partial \Gamma_k$ has finite Hausdorff measure $h_k$ 
in the gauge $\{ p_n \}$.  For each $k$ there are cutsets $\Pi_k^{(j)}$
of $\Gamma_k$ tending to infinity such that 
$$\sum_{\sigma \in \Pi_k^{(j)}} p (|\sigma|) \rightarrow h_k$$
and therefore the expected number of surviving rays 
of $\Gamma_k$ is at most $h_k$.  Since the number of surviving
rays of $\Gamma_k$ is either 0 or infinite, we conclude it is
almost surely 0.
$\Cox$

The remainder of this section makes progress
towards Conjecture~\ref{conjcapacity} by proving that positive capacity
is necessary for $\P (\rho \leftrightarrow \partial \Gamma) > 0$ 
in some useful cases.  We do this by
comparing different target percolations, varying
either $\Gamma$ or the set $B$.  The notation $p(n) = \P (\rho
\leftrightarrow \sigma)$ for $\sigma \in \Gamma_n$ is written 
$p(B;n)$ when we want to emphasize dependence on $B$;
similarly, $\P (B; \cdots)$ reflects dependence on $B$.  It is assumed
that the common distribution of the $X(\sigma)$ defining the target
percolation never change, since this may always be accomplished through
a measure-theoretic isomorphism.  It will be seen below (and this 
makes the comparison theorems useful) that spherically symmetric 
trees are much easier to handle than general trees.  This is partly
because $\P(\rho \leftrightarrow \partial \Gamma)$ may be calculated
recursively by conditioning.  In particular, if $f(n)$
is the growth function for $\Gamma$ (i.e. each $\sigma \in \Gamma_{n-1}$
has $f(n)$ neighbors in $\Gamma_n$), $\Gamma (\sigma)$ is the 
subtree of $\Gamma$ rooted at a vertex $\sigma \in \Gamma_1$,
and $B/x \subseteq \R^{N-1} = \{ y \in \R^{N-1} : (x , y_1 , y_2 , 
\ldots, y_{N-1}) \in B \}$ is the cross-section of $B$ at $x$,
then conditioning on $X(\sigma)$ for $\sigma \in \Gamma_1$ gives
\begin{equation} \label{eqrecursive}
\P (B ; \rho \notlr \partial \Gamma) = 
   \left [ \E \P ( B / X(\sigma) ; \sigma \notlr \partial 
   \Gamma (\sigma)) \right ]^{f(1)} .
\end{equation}
Notice that this recursion makes sense if $f(1)$ is any positive real,
not necessarily an integer.  As a notational convenience, we define
$\P (\rho \leftrightarrow \partial \Gamma)$ for {\em virtual} 
spherically symmetric trees with positive real growth functions, $f$,
by~(\ref{eqrecursive}) for trees of finite height and passage to the limit
for infinite trees.  Specifically, if $B$ is any target set and $\Gamma$
is any tree, we let $f(n) = |\Gamma_n| / |\Gamma_{n-1}|$ and define the symbol 
$\S(\Gamma)$ to stand for the spherical symmetrization of $\Gamma$ in the
sense that the expression $\P (B ; \rho \notlr \partial \S(\Gamma))$
is defined to stand for the value of the function $\Psi (B ; f(1) , \ldots ,
f(N))$, where $\Psi$ is defined by the following recursion in which
$X$ is a random variable with the common distribution of the $X (\sigma)$: 
\begin{quote}
$\Psi (B;a) = \P (X \notin B)^a$ if $B \subseteq \R$ and $a > 0$;

$\Psi (B ; a_1 , \ldots , a_n) = \left [ \E \Psi (B/X ; a_2 , \ldots ,
   a_n) \right ]^{a_1}$ ~~if $B \subseteq \R^n$ and $a \in (\R^n)^+$.
\end{quote}

\begin{th} \label{th4.2}
Let $\Gamma$ be any tree of height $N \leq \infty$, let $B \subset \R^N$
be any target set, and let $\S (\Gamma)$ be the (virtual) spherically
symmetric tree with $f(n) = |\Gamma_n| / |\Gamma_{n-1}|$.
Let $\S (B)$ be a Cartesian product target set $\{ y \in \R^N : y_i
\leq b_i \mbox{ for all finite } i \leq N \}$ with $b_i$ chosen
so that 
$$\prod_{i=1}^n \P (X(\sigma) \leq b_i) = p(B;n) .$$
Then
\begin{eqnarray} 
\P (B ; \rho \leftrightarrow \partial \Gamma) & \leq & \P (B ;
   \rho \leftrightarrow \partial \S (\Gamma))  \label{eqcomp1} \\[2ex]
& \leq & \P (\S (B) ; \rho \leftrightarrow \partial \S (\Gamma)) 
   \label{eqcomp2} \\[2ex]
& \leq & 2 \left [ p(B;N)^{-1} |\Gamma_N|^{-1} + \sum_{k=0}^{N-1} p(B;k)^{-1} 
   (|\Gamma_k|^{-1} - |\Gamma_{k+1}|^{-1}) \right ]^{-1} \label{eqcomp3} 
\end{eqnarray}
where the $p(B;N)^{-1} |\Gamma_N|^{-1}$ term appears only if $N < \infty$.
If $\S(\Gamma)$ exists as a tree, i.e.\ $f(n)$ is an integer
for all $n$, then the
expression~(\ref{eqcomp3}) is precisely $2 \Cap_p (\Gamma)$.  
\end{th}

When $\Gamma$ is spherically symmetric, $\S(\Gamma) = \Gamma$ and we get:
\begin{cor} \label{cor4.25}
If $\Gamma$ is spherically symmetric then $\Cap_p (\Gamma) > 0$ is necessary
for $\P (\rho \leftrightarrow \partial \Gamma) > 0$ in any target percolation.
$\Cox$
\end{cor}

\noindent{Remark and counterexample:}  The inequality $\P (B ; \rho 
\leftrightarrow \partial \Gamma) \leq \P (\S (B) ; \rho 
\leftrightarrow \partial \Gamma)$ holds for spherically symmetric trees
by~(\ref{eqcomp2}) and also for certain types of target sets $B$ (see 
Theorem~\ref{th4.4} below) but fails in general.  Whenever this 
inequality holds, Lyons' (1992) result that $\Cap_p (\Gamma) > 0$ is 
necessary for $\P (\rho \leftrightarrow \partial \Gamma) > 0$
in Bernoulli percolation implies Conjecture~\ref{conjcapacity} for that
case, since $\S (B)$ defines a Bernoulli percolation.  A counterexample
to the general inequality is the tree: \\
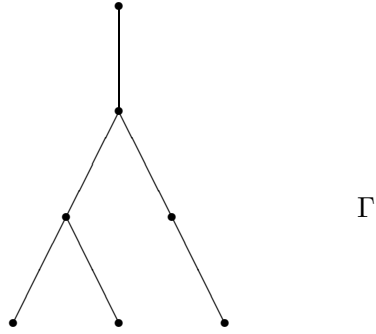
\begin{picture}(260,150)
\put(110,10){\circle*{3}}
\put(150,10){\circle*{3}}
\put(190,10){\circle*{3}}
\put(130,50){\circle*{3}}
\put(170,50){\circle*{3}}
\put(150,90){\circle*{3}}
\put(150,130){\circle*{3}}
\put(150,130){\line(0,-1){40}}
\put(150,90){\line(1,-2){20}}
\put(150,90){\line(-1,-2){20}}
\put(130,50){\line(-1,-2){20}}
\put(130,50){\line(1,-2){20}}
\put(170,50){\line(1,-2){20}}
\put(240,50){$\Gamma$}
\end{picture} \\
Let $X(\sigma)$ be uniform on $[0,1]$ and define
$$ \left ( B = [0,1/2] \times [2\ee,1] \times [0,1] \right )
\cup \left ( [1/2,1] \times [0,1] \times [4\ee,1] \right ) .$$
Then $\S (B) = [0,1] \times [0,1-\ee] \times [0 , {1 - 3\ee \over 1 - \ee}]$
and 
$$ \P (\S (B) ; \rho \leftrightarrow \partial \Gamma) \leq 1 - 3 \ee^2
  < 1 - 2 \ee^2 - 32\ee^3 = \P (B ; \rho \leftrightarrow \partial
\Gamma)$$
for sufficiently small $\ee$.  

\noindent{{\bf Question:}}  Is there an infinite tree, $\Gamma$, and a
target set, $B$, for which 
$$\P (\S (B) ; \rho \leftrightarrow \partial \Gamma) = 0 < 
  \P (B ; \rho \leftrightarrow \partial \Gamma) \; ?$$

The proof of Theorem~\ref{th4.2} is based on the following convexity lemma.

\begin{lem} \label{lem4.3}
For any tree $\Gamma$ of height $n < \infty$, define the function
$h_\Gamma (z_1 , \ldots , z_n)$ for arguments $1 \geq z_1 \geq \cdots 
\geq z_n \geq 0$ by
$$h_\Gamma (z_1 , \ldots , z_n) = \P (\S (B) ; \rho \notlr \Gamma_n)$$
where $B$ defines a target percolation with $\P (B ;  \rho \leftrightarrow
\sigma) = z_{|\sigma|}$.  If $\Gamma$ is spherically symmetric then
$h_\Gamma$ is a convex function.  The same holds for virtual trees, under
the restriction that the growth function $f(n)$ is always greater than
or equal to one.
\end{lem}

{\sc Proof:}  Proceed by induction on $n$.  When $n = 1$, certainly
$h_\Gamma (z_1) = (1-z_1)^{|\Gamma_1|}$ is convex.  Now assume the result 
for trees of height $n-1$ and let $\Gamma (\sigma)$ denote the subtree of 
$\Gamma$ rooted at $\sigma$: $\{ \tau : \tau \geq \sigma \}$.  Since
$\Gamma$ is spherically symmetric, all subtrees $\Gamma (\sigma)$
with $\sigma \in \Gamma_1$ are isomorphic, spherically symmetric
trees of height $n-1$.  By definition of $h_\Gamma$,
$$h_\Gamma (z_1 , \ldots , z_n) = \left [ (1 - z_1) + z_1 h_{\Gamma
   (\sigma)} \left ( {z_2 \over z_1} , \ldots , {z_n \over z_1} \right )
   \right ]^{|\Gamma_1|} $$
where $\sigma \in \Gamma_1$.  By induction, the function $h_{\Gamma
(\sigma)}$ is convex.  This implies that 
$$g (z_1 , \ldots , z_n) = z_1 h_{\Gamma (\sigma)} \left (
{z_2 \over z_1} , \ldots , {z_n \over z_1} \right ) $$
is convex: since $g$ is homogeneous of degree one, it suffices to check
convexity on the affine hyperplane $\{ z_1 = 1 \}$, where it is clear.
Adding a linear function to $g$ and taking a power of at least one
preserves convexity, so $h_\Gamma$ is convex, completing the induction.
$\Cox$

{\sc Proof of Theorem}~\ref{th4.2}:
The first inequality is proved in Pemantle and Peres (1994).  For the
second inequality it clearly suffices to consider trees of finite 
height, $N$.  When $N = 1$ there is nothing to prove, so fix
$N > 1$ and assume for induction that the inequality holds for trees
of height $N-1$.  Define $\Gamma ( \sigma)$, $h_\Gamma$, $B / x$ 
and $p(B;k)$ as previously, and observe that for every $j < N$, 
$$\E \, p(B / X ; j-1) = p(B ; j)$$
where $X$ has the common distribution of the $\{X(\sigma)\}$.  Use the 
induction hypothesis and the fact that $X(\sigma)$ are independent to get
\begin{eqnarray*}
&& \P (B ; \rho \notlr \partial \Gamma) \\[2ex]
& = & \prod_{\sigma \in \Gamma_1} \left [ 1 - p(B;1) + p(B;1) \, \E
   \P (B/X ; \sigma \notlr \partial \Gamma (\sigma)) \right ] \\[2ex]
& \geq & \left [ 1 - p(B;1) + p(B;1) \, \E
   h_{\Gamma (\sigma)} \left ( {p(B/X;1) \over p(B;1)} , \ldots ,
   {p(B/X;N-1) \over p(B;1)} \right ) \right ]^{|\Gamma_1|} .
\end{eqnarray*}
Utilizing the convexity of $h_{\Gamma (\sigma)}$ and Jensen's inequality,
the last expression is at least 
\begin{eqnarray*}
&& \left [ 1 - p(B;1) + p(B;1) \, h_\Gamma \left ( {p(B;2) \over p(B;1)} ,
   \ldots , {p(B;N) \over p(B;1)} \right ) \right ]^{|\Gamma_1|} \\[2ex]
& = & h_\Gamma (p(B;1) , \ldots , p(B;N)) \\[2ex]
& = & \P (\S (B) ; \rho \notlr \partial \Gamma)
\end{eqnarray*}
completing the induction and the proof of the second inequality.  

When $\Gamma$ is spherically symmetric, Theorem~2.1 of Lyons (1992)
asserts that 
$$\P ( \rho \leftrightarrow \partial \Gamma) \leq 2 \Cap_p (\Gamma).$$  
The measure $\mu$ that minimizes $I_p (\mu)$ for spherically symmetric 
trees is easily seen to be uniform, with 
$$I_p (\mu) = \int\int p(|\xi \wedge \eta|)^{-1} d\mu^2 ; $$
summing by parts shows that RHS of~(\ref{eqcomp3}) is equal to 
$2 I_p (\mu)^{-1}$.
When $\Gamma$ is not spherically symmetric, the final inequality 
is proved by induction, as follows.  

Fix $B$ and let $p_n = \P (B; \rho \leftrightarrow \sigma) = \P (\S(B) ; 
\rho \leftrightarrow \sigma)$ for $|\sigma| = n$.  
Let $\psi_p (\lambda_1 , \ldots , \lambda_N)$ denote 
$\P (\S (B) ; \rho \leftrightarrow \partial \Gamma)$ for a 
(possibly virtual) spherically symmetric tree $\Lambda$ whose growth 
numbers $f(n)$ satisfy $\prod_{i=0}^{k-1} f(i) = \lambda_k$.  Write
$R_p ( \lambda_1 , \ldots , \lambda_N)$ for the ``electrical 
resistance'' of this tree when edges at level $i$ are assigned 
resistance $p_i^{-1} - p_{i-1}^{-1}$ and an additional unit resistor
is attached to the root.  Explicitly, define
$$R_p ( \lambda_1 , \ldots , \lambda_N) = p_N^{-1} \lambda_N^{-1} +
   \sum_{i=0}^{N-1} p_i^{-1} (\lambda_i^{-1} - \lambda_{i+1}^{-1}) ,$$
where $\lambda_0 \deq 1$, so the inequality to be proved is
\begin{equation} \label{eq19}
\psi_p (\lambda_1 , \ldots , \lambda_N) \leq 2 R_p (\lambda_1 , \ldots ,
   \lambda_N)^{-1} .
\end{equation}
Proceed by induction, the case $N = 0$ boiling down to $1 \leq 2$.
Letting $p_i' = p_i / p_1$, we have
$$\psi_p (\lambda_1 , \ldots , \lambda_N)  = 1 - \left [ 1 - p_1 \psi_{p'}
    \left ( {\lambda_2 \over \lambda_1} , \ldots , {\lambda_N \over
    \lambda_1} \right ) \right ]^{\lambda_1} $$
Using the elementary inequality 
$${1 - x^{\lambda_1} \over 1 + x^{\lambda_1}} \leq \lambda_1 {1 - x
   \over 1 + x} , $$
valid for all $x \in [0,1]$ and $\lambda_1 \geq 1$, we obtain
\begin{eqnarray*}
{\psi_p ( \lambda_1 , \ldots , \lambda_N) \over 2 - \psi_p ( \lambda_1 , 
    \ldots , \lambda_N)} & = & {1 - \left [ 1 - p_1 \psi_{p'} (
    {\lambda_2 \over \lambda_1} , \ldots , {\lambda_N \over \lambda_1} )
    \right ]^{\lambda_1} \over 1 + \left [ 1 - p_1 \psi_{p'} (
    {\lambda_2 \over \lambda_1} , \ldots , {\lambda_N \over \lambda_1} )
    \right ]^{\lambda_1} } \\[3ex]
& \leq & { \lambda_1 p_1 \psi_{p'} (
    {\lambda_2 \over \lambda_1} , \ldots , {\lambda_N \over \lambda_1} )
    \over 2 - p_1 \psi_{p'} (
    {\lambda_2 \over \lambda_1} , \ldots , {\lambda_N \over \lambda_1} )} .
\end{eqnarray*}
Applying the inductive hypothesis, this is at most
$$ { 2 p_1 \lambda_1 R_{p'} (
    {\lambda_2 \over \lambda_1} , \ldots , {\lambda_N \over \lambda_1} )^{-1}
    \over 2 - 2 p_1 R_{p'} (
    {\lambda_2 \over \lambda_1} , \ldots , {\lambda_N \over \lambda_1} )^{-1}}
= {p_1 \lambda_1 \over R_{p'} (
    {\lambda_2 \over \lambda_1} , \ldots , {\lambda_N \over \lambda_1} 
    ) - p_1 } . $$
The last expression may be simplified using 
$$R_p ( \lambda_1 , \ldots , \lambda_N) = 1 - \lambda_1^{-1} + p_1^{-1}
  \lambda_1 ^{-1} R_{p'} (\lambda_2 / \lambda_1 , \ldots , \lambda_N 
  / \lambda_1)$$
to get
$$ {\psi_p ( \lambda_1 , \ldots , \lambda_N) \over 2 - \psi_p ( \lambda_1 , 
    \ldots , \lambda_N)} \leq {p_1 \lambda_1 \over p_1 \lambda_1 
    (R_p (\lambda_1 , \ldots , \lambda_N) - 1)} = {2 
    R_p (\lambda_1 , \ldots , \lambda_N)^{-1} \over 2 - 2
    R_p (\lambda_1 , \ldots , \lambda_N)^{-1}} .$$
This proves~(\ref{eq19}) and the theorem.   $\Cox$

The last result of this section gives a condition on the target set $B$,
sufficient to imply that $\P (\rho \leftrightarrow \partial \Gamma)$
increases when $B$ is replaced by $\S(B)$.  The condition is rather 
strong, however it can be applied in two very natural cases -- see
Theorem~\ref{th5.1} below.  

\begin{th} \label{th4.4}
Let $\Gamma$ be any tree of height $N \leq \infty$, let $X(\sigma)$ be
$\iid$ real random variables, and let $B$ be any target set.  For
integers $k \leq j \leq N$ and real numbers $x_1 , \ldots , x_k$, define
$$p_j (x_1 , \ldots , x_k) = \P ((X_{k+1} , \ldots ,
   X_j) \in \pi_j (B) / x_1 \cdots x_k)  $$ 
to be the measure of the cross section at $x_1 , \ldots , x_k$ of the
projection onto the first $j$ coordinates of $B$.  Suppose that
for every fixed $x_1 , \ldots , x_{k-1}$, the matrix $M$ whose $(y,j)$-entry
is $p_j (x_1 , \ldots , x_{k-1} , y)$ is totally positive of order two
($TP_2$), i.e. $M_{xi} M_{yj} \geq M_{yi} M_{xj}$ when $k \leq i < j$
and $x < y$.  Then
\begin{equation} \label{eq20}
\P (B ; \rho \leftrightarrow \partial \Gamma) \leq
    \P (\S (B)  ; \rho \leftrightarrow \partial \Gamma) .
\end{equation}
\end{th}

We shall require the following version of Jensen's inequality.

\begin{lem} \label{po convexity}
Let $\po$ be a partial order on $\R^n$ such that for every
$w \in \R^n$ the set $\{ z : z \po w \}$ is convex.  Let
$\mu$ be a probability measure supported on a bounded subset
of $\R^n$ which is totally ordered by $\po$, and let $h : \R^n
\rightarrow \R$ be a continuous function.  If $h$ is convex on
any segment connecting two comparable points $z \po w$, then
$$ h \left ( \int x \, d\mu \right ) \leq \int h(x) \, d\mu .$$
\end{lem}

{\sc Proof:}  It is enough to prove this in the case where $\mu$
has finite support, since we may approximate any measure by 
measures supported on finite subsets and use continuity of $h$
and bounded support to pass to the limit.  Letting $z_1 \po
\cdots \po z_m$ denote the support of $\mu$ and letting $a_i
= \mu \{ z_i \}$, we proceed by induction on $m$.  If $m = 2$,
the desired inequality $h ( a_1 z_1 + a_2 z_2) \leq a_1 h(z_1)
+ a_2 h(z_2)$ is a direct consequence of the assumption on $h$.
If $m > 2$, let $\nu$ be the measure which puts mass $a_1 + a_2$
at $(a_1 z_1 + a_2 z_2)/(a_1 + a_2)$ and mass $a_i$ at $z_i$ for
each $i \geq 3$.  The support of $\nu$ is a totally ordered set
of cardinality $n-1$ so the induction hypothesis implies 
$$ h \left ( \int x \, d\mu \right ) = h \left ( \int x \, d\nu \right ) 
   \leq \int h(x) \, d\nu .$$
Applying the convexity assumption on $h$ at $z_1$ and $z_2$ then gives 
$$\int h(x) \, d\nu \leq \int h(x) \, d\mu ,$$
completing the induction.    $\Cox$

{\sc Proof of Theorem}~\ref{th4.4}:  For each $n$, let $\simp_n$ denote 
the space of points $\{ (z_1 , \ldots , z_n) \in \R^n : 1 \geq
z_1 \geq \cdots \geq z_n \geq 0\}$ and for $\zz , \ww \in \simp_n$, define
$\zz \preceq \ww$ if and only if the matrix with rows $\zz$ and $\ww$
and first column $(1,1)$ is $TP_2$ (equivalently, $z_i / w_i$ is
at most one and nonincreasing in $i$).  Define $h_\Gamma$ 
on $\bigcup \simp_n$ as in Lemma~\ref{lem4.3} so that 
$$h_\Gamma (z_1 , \ldots , z_n) = \prod_{\sigma \in \Gamma_1}
   \left [ (1 - z_1) + z_1 h_{\Gamma
   (\sigma)} \left ( {z_2 \over z_1} , \ldots , {z_n \over z_1} \right )
   \right ] . $$
In order to use Lemma~\ref{po convexity} for $h_\Gamma$, $\preceq$,
and $\simp_n$, we observe first that $\preceq$ is 
closed under convex combinations in either argument.  Observe also that
$\simp_n$ is compact and $h_\Gamma$ continuous; we now
establish by induction on $n$ that $h_\Gamma$ is convex along the
line segment joining $\zz$ and $\ww$ whenever $\zz \preceq \ww \in \simp_n$.

The initial step is immediate: $h_\Gamma (z_1) = (1-z_1)^{|\Gamma_1|}$,
which is convex for $z_1 \in [0,1]$.  When $n > 1$, observe that for
each $\sigma \in \Gamma_1$, the function $1 - z_1 + z_1 h_{\Gamma (\sigma)}
 (z_2 / z_1 , \ldots , z_n / z_1)$ is decreasing along the line segment
{}from $\zz$ to $\ww$ when $\zz \preceq \ww$.  The product of decreasing
convex functions is again convex, so it suffices to check that each
$$\phi (z_1 , \ldots z_n) \deq z_1 h_{\Gamma(\sigma)} 
   (z_2 / z_1 , \ldots , z_n / z_1) $$ 
is convex along such a line segment.  Pictorially, we must show that the
graph of $\phi$ in $\simp_n \times \R$ defined by
$\{ (z_1 , \ldots , z_{n+1}) : z_{n+1} = z_1 h_{\Gamma (\sigma)} (z_2 / z_1 ,
\ldots , z_n / z_1) \}$ lies below any chord $\overline (\zz , \phi (\zz))
(\ww , \phi (\ww))$ whenever $\zz \preceq \ww$.
Observe that the graph of $\phi$ is the cone of the set $\{ (1 , z_2 , \ldots ,
z_{n+1}) : z_{n+1} = h_{\Gamma (\sigma)} (z_2 , \ldots , z_n) \}$ with
the origin.  In other words, viewing $\simp_{n-1}$ as embedded in 
$\simp_n$ by $(z_2 , \ldots , z_n) \mapsto (1 , z_2 , \ldots , z_n)$,
the graph of $\phi$ is the cone of the graph of the $n-1$-argument function
$h_{\Gamma (\sigma)}$.  To check that the chord of the graph of $\phi$
between $\zz$ and $\ww$ lies above the graph, it then suffices to
see that the chord of the graph of $h_{\Gamma (\sigma)}$ between
$(z_2 / z_1 , \ldots , z_n / z_1)$ and $(w_2 / w_1 , \ldots , w_n / w_1)$
lies above the graph.  But $\zz \preceq \ww \in \simp_n$ implies
$(z_2 / z_1 , \ldots , z_n / z_1) \preceq (w_2 / w_1 , \ldots , w_n / w_1)
\in \simp_{n-1}$, so this follows from the induction hypothesis.

We now prove the theorem for trees of finite height, the infinite
case following from writing $\P (\rho \leftrightarrow
\partial \Gamma)$ as the decreasing limit of $\P (\rho \leftrightarrow
\Gamma_n)$.  Let $N < \infty$ and proceed by induction on $N$, the
case $N=1$ being trivial since $B = \S(B)$.  Assume therefore that
$N > 1$ and that the theorem is true for smaller values of $N$.

The induction is then completed by justifying the following chain 
of identities and inequalities.  By conditioning on the independent
random variables $\{ X(\sigma) : \sigma \in \Gamma_1 \}$ and using 
the induction hypothesis we get:
\begin{eqnarray} 
\P (B ; \rho \notlr \partial \Gamma) & = & \prod_{\sigma \in \Gamma_1}
   \E \P (B / X ; \sigma \notlr \partial \Gamma (\sigma)) 
   \label{new 4.7} \\[2ex]
& \geq & \prod_{\sigma \in \Gamma_1} \E \P (\S (B) / X ; \sigma 
   \notlr \partial \Gamma (\sigma))  . \nonumber
\end{eqnarray}
Recalling the definition of $h_{\Gamma (\sigma)}$ and $p_j (x)$, this
is equal to
\begin{eqnarray}
&& \prod_{\sigma \in \Gamma_1} \E \left [ h_{\Gamma (\sigma)} (p_2 (X) , 
   \ldots , p_N (X)) \right ] \nonumber \\[2ex]
& = & \prod_{\sigma \in \Gamma_1} \left [ (1 - p_1) + p_1 \int 
   h_{\Gamma (\sigma)} (p_2 (x) , \ldots , p_N (x)) \, d\mu (x)
   \right ] \label{new 4.10} ,
\end{eqnarray}
where $\mu$ is the conditional distribution of $X$ given $X \in \pi_1
(B)$ and $p_1 = \P (X \in \pi_1 (B))$.  

Now observe that the vectors $(p_2 (x) , \ldots , p_N (x))$ with $x
\in \pi_1 (B)$ are totally ordered by $\preceq$ according to the
$k=1$ case of the $TP_2$ assumption of the theorem.  Since $\int
p_j (x) \, d\mu (x) = p_j / p_1$ for $2 \leq j \leq N$, 
Lemma~\ref{po convexity} applied to $h_{\Gamma (\sigma)}$ and
$\preceq$ on the set $\simp_{N-1}$ shows that~(\ref{new 4.10})
is at least
\begin{eqnarray}
&& \prod_{\sigma \in \Gamma_1} \left [ (1 - p_1) + p_1 h_{\Gamma (\sigma)}
   \left ( {p_2 \over p_1} , \ldots , {p_N \over p_1} \right ) \right ]
   \nonumber \\[2ex]
& = & h_{\Gamma} ( p_1 , \ldots , p_N) \nonumber \\[2ex]
& = & \P (\S (B) ; \rho \notlr \partial \Gamma )  . \label{new 4.13} 
\end{eqnarray}
Comparing~(\ref{new 4.7}) to~(\ref{new 4.13}) we see that the theorem
is established.   $\Cox$

\section{Positive rays for tree-indexed random walks}

\setcounter{equation}{0}

This section puts together results from the previous
two sections in order to prove Theorem~\ref{th2.2}
and a few related corollaries and examples.

{\sc Proof of Theorem}~\ref{th2.2}: $(i)$  We are given a tree
$\Gamma$ with positive capacity in gauge $\phi (n) = n^{-1/2}$
and $\iid$ real random variables $\{ X(\sigma) \}$ with mean zero and
finite variance.  Consider the target percolation with target set
$B^{(0)} = \{ \xx \in \R^\infty : \sum_{i=1}^n x_i \geq 0 \mbox{ for all }
n \}$.  By Proposition~\ref{prFeller},
$$c_1' |\sigma|^{-1/2} \leq \P (\rho \leftrightarrow \sigma) 
   \leq c_1'' |\sigma|^{-1/2} .$$
Now we verify that this percolation is quasi-Bernoulli, that is to say,
\begin{equation} \label{eq21}
\P ( \rho \leftrightarrow \sigma \com \tau \| \rho \leftrightarrow
\sigma \wedge \tau) \leq c \, { |\sigma \wedge \tau| \over |\sigma|^{1/2}
   |\tau|^{1/2}} 
\end{equation}
for $\sigma , \tau \in \Gamma$.  Assume without loss of generality that
$|\sigma \wedge \tau| \leq (1/2) \min (|\sigma| , |\tau|)$, since
otherwise the claim is immediate.  The LHS of~(\ref{eq21}) is equal to
\begin{equation} \label{eq 940102b}
\begin{array}{l}
\int_0^\infty  \P (S(\sigma \wedge \tau) \in dy \| \rho \leftrightarrow
  \sigma \wedge \tau) \cdot \\
  ~~~~\P (\rho \leftrightarrow \sigma \| \rho \leftrightarrow \sigma \wedge 
     \tau, S(\sigma \wedge \tau) = y)  \cdot 
  \P (\rho \leftrightarrow \tau \| \rho \leftrightarrow \sigma \wedge 
     \tau, S(\sigma \wedge \tau) = y)  \, . \end{array}
\end{equation}
Recalling the definition of $T_y$ as the first time a trajectory
is less than $-y$, we may write
\begin{eqnarray*}
\P (\rho \leftrightarrow \sigma \| \rho \leftrightarrow \sigma \wedge 
   \tau, S(\sigma \wedge \tau) = y) & = & \\
\P ( T_y \geq |\sigma| - |\sigma \wedge \tau|) & \leq &
   4 c_2 {y \over |\sigma|^{1/2}} 
\end{eqnarray*}
by part~$(i)$ of Lemma~\ref{lem ii-iv}.  A similar bound holds for the
last factor in~(\ref{eq 940102b}).  Now use the the second part 
of Lemma~\ref{lem ii-iv} to show that~(\ref{eq 940102b}) is at most
\begin{eqnarray*}
&& 4 c_2^2 \int_0^\infty  \P (S(\sigma \wedge \tau) \in dy \| \rho
    \leftrightarrow \sigma \wedge \tau) {y^2 \over |\sigma|^{1/2}
    |\tau|^{1/2}} \\[2ex]
& = & {4 c_2^2 \over |\sigma|^{1/2} |\tau|^{1/2}} \E (S(\sigma
   \wedge \tau)^2 \| \rho \leftrightarrow \sigma \wedge \tau) \\[2ex]
& \leq & 4 c_2^2 c_3 \, {|\sigma \wedge \tau| \over |\sigma|^{1/2}
   |\tau|^{1/2}} ,
\end{eqnarray*}
verifying the claim.

Putting $|\sigma| = |\tau| = n$ and $|\sigma \wedge \tau| = k$, 
it immediately follows that $\P (\rho \leftrightarrow \sigma \com \tau)
\leq c \sqrt{k} / n$, and the second moment method (Lemma~\ref{lem4.1}
part~$(ii)$) implies that with positive probability a ray exists
along which $S(\sigma)$ remains nonnegative.  To obtain the full
assertion of the theorem, let $f(n)$ be any increasing sequence
satisfying $\sum n^{-3/2} f(n) < \infty$ and define a new target
percolation with target set 
$$B^{(f)} = \{ \xx \in \R^\infty : \sum_{i=1}^n x_i \geq f(n)
   \one_{n \geq n_f} \mbox{ for all } n \geq 1 \} $$
where $n_f$ is as in Theorem~\ref{th3.1}, part~$(i)$.  The
conclusion of Theorem~\ref{th3.1}, part~$(i)$, shows that
$\P (B^{(f)} ; \rho \leftrightarrow \sigma)$ is of order
$|\sigma|^{-1/2}$; since $\P (B^{(f)} ; \rho \leftrightarrow
\sigma \com \tau) \leq \P (B^{(0)} ; \rho \leftrightarrow
\sigma \com \tau)$, the new percolation $B^{(f)}$ is still quasi-Bernoulli.
Thus with positive probability, $\Gamma$ contains a ray along which
$\S(\sigma) \geq f(|\sigma|)$ for sufficiently large $|\sigma|$.
To see that this event actually has probability one, observe that
it contains the tail event
$$\{ \exists \xi \in \partial \Gamma \; \exists C > 0 \; \forall
   \sigma \in \xi , \; S(\sigma) \geq f(|\sigma|) + |\sigma|^{1/4}
   - C \} , $$
which has positive probability by the preceding argument, hence 
probability one.

$(ii)$ Assume that $\Gamma$ has $\sigma$-finite Hausdorff measure
in gauge $\phi (n) = n^{-1/2}$.  For any nondecreasing function, $f$,
consider the random subgraph $W_{-f} = \{ \sigma \in \Gamma : S(\sigma)
\geq -f (|\sigma|)\}$.  By the summability assumption on $f$ and by
Theorem~\ref{th3.1} part II, there is a
constant $c$ for which  
$$p(|\sigma|) = \P (\rho \leftrightarrow \sigma) \leq c |\sigma|^{-1/2} .$$
The sharpened first moment method (Lemma~\ref{sharpened}) implies that
$W_{-f}$ almost surely fails to contain a ray of $\Gamma$.  This easily
implies the stronger statement that the subgraph $W_{-f}$ has no infinite
components, almost surely.

$(iii)$  Define $W_{-f}$ and $p(|\sigma|)$ as above.  From
Theorem~\ref{th4.2} we get
\begin{equation} \label{eq23}
\P (\rho \leftrightarrow \partial \Gamma) \leq 2 \left [ \sum_{k=0}^\infty
    p(k)^{-1} \left ( |\Gamma_k|^{-1} - |\Gamma_{k+1}|^{-1} \right )
    \right ]^{-1} .
\end{equation}
Since $p(k) \leq c k^{-1/2}$, summation by parts shows that
\begin{eqnarray*}
&& \sum_{k=0}^\infty p(k)^{-1} \left ( |\Gamma_k|^{-1} - |\Gamma_{k+1}|^{-1} 
    \right ) \\[2ex]
& \geq & c^{-1} \sum_{k=0}^\infty k^{1/2} \left ( |\Gamma_k|^{-1} -
    |\Gamma_{k+1}|^{-1} \right ) \\[2ex]
& = & c^{-1} \sum_{k=1}^\infty [k^{1/2} - (k-1)^{1/2}] |\Gamma_k|^{-1}
    \\[2ex]
& \geq & (2c)^{-1}  \sum_{k=1}^\infty k^{-1/2} |\Gamma_k|^{-1} ,
\end{eqnarray*}
so if the last sum is infinite then the RHS of~(\ref{eq23}) is
zero, completing the proof.    $\Cox$

For certain special distributions of the step sizes
$\{ X(\sigma)\}$, the class of trees for which some ray stays
positive with positive probability may be sharply delineated.

\begin{th} \label{th5.1}
Let $\Gamma$ be any infinite tree and let the $\iid$ random variables
$\{X(\sigma)\}$ have common distribution $F_1$ or $F_2$, where
$F_1$ is a standard normal and $F_2$ is the distribution putting
probability $1/2$ each on $\pm 1$.  Then the probability that 
$S(\sigma) \geq 0$ along some ray of $\Gamma$ is nonzero if and
only if $\Gamma$ has positive capacity in gauge $\phi (n) = n^{-1/2}$.
\end{th}

{\bf Remark:}  The usual variants also follow.  When $\Gamma$ has
positive capacity in gauge $\phi$, the probability is one that some
ray of $\Gamma$ has $S(\sigma) < 0$ finitely often.  This is equivalent
to finding, almost surely, a ray for which $S(\sigma) \geq f(|\sigma|)$
all but finitely often, for any monotone $f$ satisfying $\sum
n^{-3/2} |f(n)| < \infty$.

{\sc Proof:}  Let $B \subseteq \R^\infty$ be the target set $\{
\xx \in \R^\infty : \sum_{i=1}^n x_i \geq 0 \mbox{ for all } n \}$.  
One half of the theorem, namely that positive capacity implies
$\P (B ; \rho \leftrightarrow \partial \Gamma) > 0$, follows
immediately from part~$(i)$ of Theorem~\ref{th2.2}.  For the 
other half, observe that zero capacity implies
$\P (\S (B) ; \rho \leftrightarrow \partial \Gamma) = 0$,
since $\S(B)$ is Bernoulli with the same values of $p(n)$,
and $\Cap_p (\Gamma) > 0$ is known to be necessary and sufficient
for percolation (c.f. remarks after Conjecture~\ref{conjcapacity}).
The present theorem then follows from~(\ref{eq20}) once the
conditions of Theorem~\ref{th4.4} are verified.  It suffices
to establish $M_{yj} / M_{xj} > M_{yi} / M_{xi}$ for $y > x$,
in the case where $j = i+1$.

To verify this, pick any $j,k$ and any $x_1 , \ldots , x_k \geq 0$
and observe that $p_j (x_1 , \ldots , x_k) = \P (x_1 + \cdots + x_k + 
S_i \geq 0 \mbox{ for all } i \leq j)$, where $\{ S_i \}$ is a random 
walk with step sizes distributed as the $\{ X(\sigma)\}$.  
It suffices then to show that for $0 \leq x < y$,
$$\P (y + S_{j+1} \geq 0 \| y + S_i \geq 0 : i \leq j) \geq
  \P (x + S_{j+1} \geq 0 \| x + S_i \geq 0 : i \leq j) .$$

To see this in the case of $F_1$, use induction on $j$.  For $j = 0$
one pointmass obviously dominates the other.  Assuming it now
for $j-1$, write 
$$\P (y + S_{j+1} \geq 0 \| y + S_i \geq 0 : i \leq j) =
   \int d\nu_y (z) \P (z + S_j \geq 0 \| z + S_i \geq 0 : i \leq j-1) ,$$
where $\nu_y (z)$ is the conditional measure of $y + S_1$ given 
$S_i \geq -y : i \leq j$.  The Radon-Nikodym derivative
$d\nu_y / d\nu_x$ at $z \geq y$ is $dF_1 (z-y) / dF_1 (z-x)$ times
a normalizing constant.  This is an increasing function of $z$,
by the increasing likelihood property of the normal distribution.
Thus $\nu_y$ stochastically dominates $\nu_x$.  By induction, 
the integrand is increasing in $z$, which, together with the
stochastic domination, establishes the inequality.  The same
argument works for $F_2$, noting that $y - x$ is always an even
integer.    $\Cox$

\begin{cor} \label{bouncing rays}
Suppose the edges of an infinite tree $\Gamma$ are labeled by $\iid$,
mean zero, finite variance random variables $\{ X(\sigma)\}$ with 
partial sums $\{ S(\sigma)\}$.  Assume that
$$ \begin{array}{l} \mbox{either } \Gamma \mbox{ is spherically
 symmetric} \\ \mbox{or } \\ \{ X(\sigma)\} \mbox{ are normal or
take values } \pm 1 \end{array}  \hspace{2in} (*)$$
If there is almost surely a ray along which $\inf S(\sigma) > - \infty$,
then there is almost surely a ray along which $\lim S(\sigma) = + \infty$.
\end{cor}

{\sc Proof:}  Both are equivalent to $\Gamma$ having positive capacity in
gauge $n^{-1/2}$.   $\Cox$

{\bf Problem:}  Remove the assumption (*).  

If the moment generating function
of $X$ fails to exist in a neighborhood of zero, it is
possible that $\E X < 0$ but still some trees of polynomial
growth have rays along which $S(\sigma) \rightarrow \infty$.  The
critical growth exponent need not be $1/2$ in this case.  We conclude
this section with such an example.

Suppose that the common distribution of the $X(\sigma)$ is a symmetric,
stable random variable with index $\alpha \in (1,2)$.  Fix $c > 0$
and consider the target set 
$$B = \{ \xx \in \R^\infty : \sum_{i=1}^n x_i > cn \mbox{ for all } n \}$$
and 
$$B' = \{ \xx \in \R^\infty : cn^2 > \sum_{i=1}^n x_i > cn 
    \mbox{ for all } n \} . $$
The following estimates may be proved.

\begin{pr} \label{pr stables}
For $\sigma , \tau \in \Gamma_n$, let $k = |\sigma \wedge \tau|$.  Then
\begin{equation} \label{eqst1}
\P (B' ; \rho \leftrightarrow \sigma) \leq \P (B ; \rho \rightarrow \sigma) 
    \leq c_1 n^{-\alpha} 
\end{equation}
\begin{equation} \label{eqst2}
\P (B' ; \rho \leftrightarrow \sigma \com \tau) / \P (B' ; \rho \rightarrow
   \sigma)^2 \leq c_2 (\ee) k^{1 + 4\alpha + \ee}
\end{equation}
for any $\ee > 0$ and some constants $c_i > 0$.
\end{pr}   $\Cox$

Suppose that $\Gamma$ is a spherically symmetric tree with growth rate
$|\Gamma_n| \approx n^{\beta}$.  Plugging~(\ref{eqst1}) into 
Lemma~\ref{lem4.1}, we see that $\P (B' ; \rho \leftrightarrow \partial 
\Gamma)$ is zero when $\beta < \alpha$, while plugging in~(\ref{eqst2})
shows that $\P (B' ; \rho \leftrightarrow \partial 
\Gamma)$ is positive when $\beta > (1 + 4 \alpha)$. 
If one then considers the tree-indexed random walk whose increments
are distributed as $X - c$, one sees that $\beta < \alpha$ implies
that with probability one $S(\sigma) < 0$ infinitely often on every ray,
whereas $\beta > 1 + 4 \alpha$ implies that with probability one
$S(\sigma) \rightarrow \infty$ with at least linear rate along some ray.
For $\beta > 1 + 4 \alpha$, RWRE with this distribution of $X$
is therefore transient even though $\E X < 0$.

Defining the {\em sustainable speed} of a tree-indexed random walk 
to be the almost surely constant value 
$$\sup_{\xi} \liminf_{\sigma \in \xi} {S(\sigma) \over |\sigma|} ,$$
Lyons and Pemantle (1992) have shown that the Hausdorff dimension of 
$\Gamma$ and the distribution of $X$ together determine the sustainable 
speed of 
the tree-indexed random walk, as long as $X$ has a moment generating
function in a neighborhood of zero.  When the increments are symmetric
stable random variables, the moment hypothesis is violated, and the 
analysis above shows that the sustainable speed of the tree-indexed 
random walk can be different
for different polynomially growing trees of Hausdorff dimension zero.

\section{Critical RWRE: proofs}

\setcounter{equation}{0}

The following easy lemma will be useful.

\begin{lem} \label{bottleneck}
If $\Gamma$ is any tree with conductances $C(\sigma)$, let
$U(\sigma) = \min_{\rho < \tau \leq \sigma} C(\tau)$.  
Then the conductance from $\rho$ to a cutset $\Pi$ is
at most 
\begin{equation} \label{qottle}
\sum_{\sigma \in \Pi} U(\sigma) .
\end{equation}
\end{lem}

{\sc Proof:}  For each $\sigma \in \Pi$, let $\gamma (\sigma)$
be the sequence of conductances on the path from $\rho$ to $\sigma$, 
and let $\Gamma'$ be a tree consisting of disjoint paths for each $\sigma 
\in \Pi$, each path having conductances $\gamma (\sigma)$.  $\Gamma$
is a contraction of $\Gamma'$, so by Rayleigh's monotonicity law
(Doyle and Snell 1984), the
conductance to $\Pi$ in $\Gamma$ is less than or equal to the conductance
of $\Gamma'$, which is the sum over $\sigma \in \Pi$ of conductances
bounded above by $U(\sigma)$.     $\Cox$

{\sc Proof of Theorem}~\ref{th2.1}:
$(i)$  This is almost immediate from Theorem~\ref{th2.2}, which was proved
in the previous section.  Since $\Gamma$ has positive capacity in gauge
$n^{-1/2}$, that theorem guarantees the almost sure existence of a
ray $\xi$ along which the partial sums $S(\sigma) = \sum_{\rho <
\tau \leq \sigma} X(\tau)$ satisfy $S(\sigma) \geq 2 \log |\sigma|$ when
$\sigma$ is sufficiently large.  The total resistance along
$\xi$ is then $\sum_{\sigma \in \xi} e^{-S(\sigma)} < \infty$, so
the resistance of the entire tree is finite and the RWRE transient.

$(ii)$  First we calculate a bound for the expected conductance along the
path from $\rho$ to $\sigma$ in terms of $|\sigma|$, then use the Hausdorff
measure assumption to bound the net conductance of the tree from $\rho$ to
$\partial \Gamma$.  In this calculation it is convenient to attach an
auxiliary unit resistor, thought of as an edge $\rho' \rho$ comprising the
$-1$ level of $\Gamma$.  After this addition, the minimal conductance in
any edge in the path from $\rho'$ to $\sigma$ is $e^{m(\sigma)}$, where
$m(\sigma) = \min \{ S(\tau) : \rho \leq \tau \leq \sigma \}$ is
nonpositive since $S(\rho) = 0$.  Applying the first part of
Lemma~\ref{lem ii-iv} gives 
$$\P (m(\sigma) \geq -y) \leq c |\sigma|^{-1/2} \max \{ y,1 \} $$
for some constant $c$.  Plug this into the identity
$$\E e^{m(\sigma)} = \int_0^1 \P (e^{m(\sigma)} \geq u) \, du ,$$
changing variables to $y = \log (1/u)$ to get
\begin{eqnarray*}
&& \E e^{m(\sigma)} = \int_0^\infty \P (m(\sigma) \leq -y) e^{-y} 
   \, dy \\[2ex]
& \leq & c |\sigma|^{-1/2} \left ( 1 + \int_1^\infty y e^{-y} \, dy 
   \right ) \\[2ex]
& \leq & 2 c |\sigma|^{-1/2} .
\end{eqnarray*}
The Hausdorff measure assumption implies the for any $\ee > 0$ there is
a cutset $\Pi (\ee)$ for which $\sum_{\sigma \in \Pi (\ee)} |\sigma|^{-1/2}
< \ee$, hence by Lemma~\ref{bottleneck}
the expected conductance from $\rho'$ to $\Pi (\ee)$ is at most $2c\ee$.
Thus the net conductance from $\rho'$ to $\partial \Gamma$ vanishes almost
surely, so the RWRE is almost surely recurrent.

$(iii)$  Set $f(n) = \log |\Gamma_n| + 2 \log (n+1)$.  The assumptions 
in~$(iii)$ imply that $f$ is increasing and $\sum n^{-3/2} f(n) < \infty$,
so Theorem~\ref{th2.2} implies that with probability one, every ray
of $\Gamma$ has $S(\sigma) < -f(|\sigma|)$ infinitely often.  Thus
for every $N \geq 1$ there exists almost surely a random cutset $\Pi$
such that for every $\sigma \in \Pi$, $|\sigma| > N$ and $S(\sigma)
< -f(|\sigma|)$.  The net conductance from $\rho$ to this $\Pi$ is 
at most 
\begin{eqnarray*}
\sum_{\sigma \in \Pi} e^{S(\sigma)} & \leq & \sum_{\sigma \in \Pi}
   e^{-f(|\sigma|)}  \\[2ex]
& \leq & \sum_{n=N}^\infty |\Gamma_n| e^{-f(n)} \\[2ex]
& \leq & \sum_{n=N}^\infty n^{-2} .
\end{eqnarray*}
Taking $N$ large shows that the net conductance to $\partial \Gamma$
vanishes almost surely.   $\Cox$

\section{Reinforced random walk}
 Reinforced RW is a  process
 introduced by Coppersmith and Diaconis (unpublished) to model a tendency
of the random walker to revisit familiar territory.
The variant of this process analysed in Pemantle (1988) has
an inherent bias toward the root, i.e. a positive backward push;
 here we consider the following
``unbiased'' variant which fits into the general framework of reinforcement
described in Davis (1990), and may be analysed using the tools developed
in the previous sections of this paper. 

Let $\Gamma$ be an infinite rooted tree, with (dynamically changing) 
positive weights $\, w_n(e) $ for $ n\geq 0$ attached to each edge $e$.
 At time zero, all weights are set to one:
$w_0(e)=1$ for all $e$. Let $Y_0$ be the root of $\Gamma$, and for
every $n \geq 0$, given $Y_1, \ldots Y_n$ let $Y_{n+1}$ be a randomly chosen vertex
adjacent to $Y_n$, so that each edge $e$ emanating from $Y_n$
has conditional probability proportional to $w_n(e)$ to be the  edge  
connecting $Y_n$ to $Y_{n+1}$. Each time an edge is traversed
{\em  back and forth}, its weight is increased by 1, i.e. $w_k(e)-1$ is
the number of ``return trips taken on $e$ by time $k$''.
Call the resulting process $\{Y_n\}$ an ``unbiased reinforced random walk''.

\begin{th} \label{th7.1} 
\begin{description}
\item{$(i)$}  If $\Gamma$ has positive capacity in gauge $\phi(n) = n^{-1/2}$,
then the resulting reinforced RW is transient, i.e.,
 $\P(Y_n=Y_0 \mbox{\, infinitely often})=0$.
\item{$(ii)$}  If $\Gamma$ has zero Hausdorff measure in the
same gauge, then the reinforced RW  is recurrent, i.e.,
  $\P(Y_n=Y_0 \mbox{\, infinitely often})=1$.
\end{description}
\end{th}

\noindent{\sc Proof:}  Fix a vertex $\sigma$ in $\Gamma$,
of degree d. As explained in Section 3 of Pemantle (1988),
 for every vertex $\sigma$ of $\Gamma$,
the sequence of edges by which the walk leaves $\sigma$ is equivalent
to ``Polya's urn'' stopped at a random time.  Initially the urn contains
d balls, one of each color. (The colors correspond to the edges emanating 
from $\sigma$.) Each time the walk leaves $\sigma$, a ball is picked at random
from the urn,
and returned to the urn along with another ball of the same color.
(This corresponds to increasing the weight of the relevant edge).
>From Section VII.4 of Feller (1966) we find that the sequence of
edges taken from $\sigma$ is a stochastically equivalent to
a mixture of sequences of i.i.d. variables, where the mixing measure is
uniform over the simplex of probability vectors of length $d$.
A standard method to generate a uniform random  vector on the 
simplex is to pick $d \,$ independent identically distributed
exponential random variables, and normalize them by their sum.
This leads to the following 
{\bf RWRE  description of the reinforced RW}: 
\newline Assign to each edge $e$ in 
$\Gamma$ two exponential random variables $U( \stackrel{\rightarrow}{e})$
and $U( \stackrel{\leftarrow}{e})$,
 one for each orientation, so that
all the assigned variables are i.i.d. .
These labels are then used to define an environment for a random walk
on $\Gamma$ such that the transition probability from a vertex $\sigma$
to a neighboring vertex $\tau$ is 
$$
q(\sigma,\tau)={U(\stackrel{\rightarrow}{\sigma \tau}) \over 
  \sum \{U( \stackrel{\rightarrow}{e}) : \;
  \stackrel{\rightarrow}{e} \mbox{\, emanates from \,} \sigma \}  } \; .
$$
Thus the log-ratios $\{X(\sigma) \}_{\sigma \in \Gamma}$ 
defined in (\ref{eq1}) are identically distributed, and any subcollection
of these variables where no two of the corresponding vertices
are siblings, is independent. Clearly the variables $\{X(\sigma) \}$
have mean 0 and finite variance.  

The proof of Theorem 2.1$(ii)$ goes over unchanged to prove part $(ii)$
of the present theorem, since in order to apply  the "first moment method",
Lemma 4.1$(i)$, it suffices that for any ray $\xi$ in $\Gamma$, 
the variables $\{X(\sigma) \}$ for $\sigma$ on $\xi$ be independent.

To prove part $(ii)$ via the second moment method, Lemma 4.1$(ii)$,
consider the percolation process defined by retaining
only vertices for which the partial sum from the root of $X(\sigma)$
is positive. It suffices to verify that this percolation is quasi-Bernoulli,
i.e. it satisfies (5.1);  this involves only a minor modification (which we omit)
of the proof of Theorem 2.2$(i)$ given in section 5.
$\Cox$

\noindent{\sc Acknowledgements:} We are indebted to the referee for a remarkably 
careful reading of the paper. 
We gratefully acknowledge Russell Lyons and the Institute for 
Iterated Dining for bringing us together.  

\renewcommand{\baselinestretch}{1.0}\large

\end{document}